\documentclass[11pt]{amsart}
\usepackage[T1]{fontenc}
\usepackage[letterpaper,textwidth=6.25in,textheight=8.7in,centering]{geometry}
\usepackage{amsmath,amssymb,amsthm,mathtools,xcolor}
\usepackage{booktabs}
\usepackage[hidelinks]{hyperref}
\hypersetup{pdftitle={The optimal constants in Bennett's bilinear inequality},pdfauthor={Daniel Nunez-Alarcon, Daniel M. Pellegrino, Anselmo Raposo Jr., Diana Marcela Serrano-Rodriguez},pdfsubject={Bennett's inequality, sharp constants and matrix signings},pdfkeywords={Bennett's inequality, Kahane--Salem--Zygmund inequality, sharp constants, matrix signings}}
\newtheorem{theorem}{Theorem}[section]
\newtheorem{lemma}[theorem]{Lemma}
\newtheorem{corollary}[theorem]{Corollary}
\newtheorem{proposition}[theorem]{Proposition}

\newcommand{\norm}[1]{\left\|#1\right\|}

\numberwithin{equation}{section}
\title[Optimal constants in Bennett's bilinear inequality]{The optimal constants in Bennett's bilinear inequality}
\date{}
\author[D. N\'u\~nez-Alarc\'on]{Daniel N\'u\~nez-Alarc\'on}
\address{Departamento de Matem\'aticas, Facultad de Ciencias, Universidad Nacional de Colombia, Sede Bogot\'a, Bogot\'a, Colombia}
\email{dnuneza@unal.edu.co}
\author[D. M. Pellegrino]{Daniel M. Pellegrino}
\address{Departamento de Matem\'atica, Universidade Federal da Para\'iba, 58051-900 Jo\~ao Pessoa, PB, Brazil}
\email{daniel.pellegrino@academico.ufpb.br}
\author[A. Raposo Jr.]{Anselmo Raposo Jr.}
\address{Departamento de Matem\'atica, Universidade Federal do Maranh\~ao, 65085-580 S\~ao Lu\'is, MA, Brazil}
\email{anselmo.junior@ufma.br}
\author[D. M. Serrano-Rodr\'iguez]{Diana Marcela Serrano-Rodr\'iguez}
\address{Departamento de Matem\'aticas, Facultad de Ciencias, Universidad Nacional de Colombia, Sede Bogot\'a, Bogot\'a, Colombia}
\email{diserranor@unal.edu.co}

\begin{document}
\begin{abstract}
\begingroup
We determine the optimal constants in Bennett's real bilinear inequality for $p,q\in[2,\infty]$ in terms of the norms of five fixed sign matrices and their transposes, with at most five rows and five columns. In particular, the uniform bound $\sqrt{8/5}$ of Pellegrino and Raposo is sharp. We also compute the constants on the boundary with an infinite exponent and identify the unique transition between the two extremal sizes.
\endgroup
\end{abstract}
\maketitle
\section{Introduction}
\begingroup
We determine the sharp dimension-independent constants in Bennett's real bilinear inequality \cite{Bennett} for $p,q\in[2,\infty]$ by reducing the problem to five explicit matrices. In particular, the uniform bound $\sqrt{8/5}$ obtained in \cite{PR2} is sharp.
\endgroup The problem of choosing signs so that a polynomial with prescribed coefficient moduli has small supremum norm goes back to Salem and Zygmund's work on random trigonometric series \cite{SZ} and to Kahane's subsequent developments \cite{Kahane}. In the multilinear setting, the Kahane--Salem--Zygmund inequality provides forms with coefficients of modulus one whose norms have the least possible order of growth in the dimensions. Such forms play an important role in the study of coefficient inequalities and tensor products. Related matrix estimates were obtained by Bennett, Goodman and Newman \cite{BGN}, and the connection with tensor products and operator theory can be traced to \cite{MT,Var}. Polynomial versions and asymptotic estimates for unimodular multilinear forms were further studied in \cite{AR,Boas,MS,PRS}. Vector-valued extensions with asymptotically bounded constants appear in \cite{RT}.

In his 1977 paper on Schur multipliers, Bennett \cite[Proposition 3.2]{Bennett} proved a bilinear version of the Kahane--Salem--Zygmund inequality in which the dimensions and the exponents of the two spaces need not coincide. In particular, given $p_1,p_2\in[2,\infty]$ and positive integers $n_1,n_2$, there exist signs $\varepsilon_{ij}\in\{-1,1\}$ such that the bilinear form
\[
 T:\ell_{p_1}^{n_1}\times\ell_{p_2}^{n_2}\longrightarrow\mathbb R,
 \qquad T(x,y)=\sum_{i=1}^{n_2}\sum_{j=1}^{n_1}\varepsilon_{ij}x_jy_i
\]
satisfies
\begin{equation}\label{intro:bennett}
 \norm{T}\le C_{p_1,p_2}n_1^{1/2-1/p_1}n_2^{1/2-1/p_2}
 \max\{\sqrt{n_1},\sqrt{n_2}\}.
\end{equation}
Here $1/\infty=0$, $\norm{T}=\sup_{\norm{x}_{p_1}\le1,\,\norm{y}_{p_2}\le1}|T(x,y)|$, and $C_{p_1,p_2}$ is independent of $n_1,n_2$. If $n_1=n_2=n$, the dimensional factor in \eqref{intro:bennett} reduces to $n^{3/2-1/p_1-1/p_2}$. Bennett's inequality also includes exponents below two; we consider the full range $[2,\infty]^2$.

The dimensional powers in \eqref{intro:bennett} are optimal. The study of the constants multiplying them led to deterministic constructions based on Hadamard matrices. In \cite{PR} the normalized minima were shown to be asymptotically bounded above by one as the dimensions grow. Subsequently, in \cite{PR2}, a uniform upper bound of $\sqrt{8/5}$ was obtained for the constants in \eqref{intro:bennett} in the range $[2,\infty]^2$. The optimality of this uniform bound and the dependence of the least constant on $(p_1,p_2)$ remained open. 

We prove that $\sqrt{8/5}$ is optimal and determine the least constant $C_{p_1,p_2}$ for every $(p_1,p_2)\in[2,\infty]^2$. The uniform bound is attained when $p_1=p_2=2$, already in dimension $5\times5$. At $(2,\infty)$ and $(\infty,2)$ the constant is $\sqrt{37}/5$, whereas at $(\infty,\infty)$ it is $3\sqrt3/5$ (Theorem~\ref{thm:boundary}). \begingroup
Theorem~\ref{thm:continuous} reduces the calculation of $C_{p_1,p_2}$ to the norms of two $5\times5$ matrices, two $5\times4$ matrices, one $5\times3$ matrix, and their transposes. These matrices are independent of $p_1$ and $p_2$. On the boundary with an infinite exponent the formula involves only two elementary functions of the remaining exponent.
\endgroup
 If complex coefficients of modulus one are allowed, the optimal dimension-independent constant is one for every pair of exponents \cite{PRS}: a submatrix of the (unnormalized) Fourier matrix of order $\max(n_1,n_2)$ gives the upper bound, and the case $n_1=n_2=1$ gives the lower bound.

\begingroup
\subsection*{The optimal uniform constant}
The upper estimate of Pellegrino and Raposo \cite{PR2} is already sharp at $p_1=p_2=2$. Indeed, let $A$ be any $5\times5$ sign matrix. Multiplication of individual rows by $-1$ preserves its spectral norm and allows us to assume that each row has product one. The Hamming distance between two rows is then even. Their inner product is $5-2d$, where $d\in\{0,2,4\}$ is that distance. If two rows coincide, the corresponding principal submatrix of $AA^t$ is $\left(\begin{smallmatrix}5&5\\5&5\end{smallmatrix}\right)$ and has largest eigenvalue $10$. If two rows differ in four positions, that submatrix is $\left(\begin{smallmatrix}5&-3\\-3&5\end{smallmatrix}\right)$ and has largest eigenvalue $8$. In either case the Rayleigh quotient, evaluated on a vector supported on the two rows, gives $\|A\|_{2\to2}^2\ge8$. In the remaining case every two distinct rows differ in two positions, so $AA^t=4I_5+J_5$. Its eigenvalues are $9$ and $4$ (the latter with multiplicity four). Thus in all cases $\|A\|_{2\to2}\ge\sqrt8$. Since the dimensional factor in \eqref{intro:bennett} is $\sqrt5$ at $(2,2)$, we obtain
\[
 \sqrt{8/5}\le C_{2,2}\le\sup_{p_1,p_2\in[2,\infty]}C_{p_1,p_2}\le\sqrt{8/5},
\]
where the last inequality is \cite{PR2}. This proves the optimality of the uniform constant without any classification of sign matrices.

The classical Hadamard constructions used here originate with Sylvester and Paley \cite{Paley}; see \cite{Horadam} for a systematic account. At the single endpoint $(\infty,\infty)$, minimizing the norm in a prescribed size is the matrix formulation of the Gale--Berlekamp switching game: it is the minimum, over sign matrices $A$, of $\max_{u,v\text{ signs}}v^tAu$. Brown and Spencer \cite{BrownSpencer} studied its order of growth, and Fishburn and Sloane \cite{FishburnSloane} and Carlson and Stolarski \cite{CarlsonStolarski} studied the finite square game. The estimates in \cite{PR2} apply to both Bennett's inequality and that game. Here the optimization is over all rectangular sizes and both exponents; the reduction to an exponent-independent family of fixed sign matrices is the principal result.
\endgroup
\subsection*{Reduction to finitely many sizes}
Write
\[
 d_{p_1,p_2}(n_1,n_2)=n_1^{1/2-1/p_1}n_2^{1/2-1/p_2}\sqrt{\max(n_1,n_2)}
\]
and denote by $k_{p_1,p_2}(n_1,n_2)$ the minimum of $\norm{T}/d_{p_1,p_2}(n_1,n_2)$ over the $2^{n_1n_2}$ possible choices of signs. Thus
\begin{equation}\label{intro:supmin}
 C_{p_1,p_2}=\sup_{n_1,n_2\ge1}k_{p_1,p_2}(n_1,n_2).
\end{equation}
A size $(n_1,n_2)$ cannot attain the supremum in \eqref{intro:supmin} if a single signing in that size has normalized norm strictly smaller than the normalized norm of every signing in another size.

There are only four sign vectors of length three modulo multiplication by $-1$. Among five rows, two must coincide up to sign. If $u$ is one of these rows, two coordinates of $Au$ have modulus three and each of the other three has modulus at least one. Thus $\|Au\|_1\ge9$. Testing the bilinear form on $3^{-1/p_1}u$ and on $5^{-1/p_2}$ times the coordinatewise signs of $Au$ gives
\begin{equation}\label{intro:gap}
 C_{p_1,p_2}\ge k_{p_1,p_2}(3,5)\ge\frac{9}{5\sqrt3}=\sqrt{27/25}>1.
\end{equation}
If $H$ is a Hadamard matrix of order $h$, so that $H^tH=hI_h$, then any $n_2\times n_1$ submatrix $A$ of $H$, with $h\ge N=\max(n_1,n_2)$, satisfies
\begin{equation}\label{intro:had}
 \|A\|_{p_1\to p_2'}\le \sqrt h\,n_1^{1/2-1/p_1}n_2^{1/2-1/p_2},\qquad
 k_{p_1,p_2}(n_1,n_2)\le\sqrt{h/N}.
\end{equation}
The first inequality follows by applying the inclusions of the $\ell_p$ balls into the corresponding scaled $\ell_2$ balls and the identity $\|H\|_{2\to2}=\sqrt h$. Consequently, a dimension is excluded as soon as there is a Hadamard order $h\ge N$ with
\begin{equation}\label{intro:threshold}
 \boxed{\frac hN<\frac{27}{25}.}
\end{equation}
The Hadamard orders in Section~\ref{sec:had} cover the interval $[104,208]$ with consecutive ratios at most $14/13<27/25$. Tensor products with Sylvester matrices of orders $2^j$ cover every dyadic dilation of that interval. Together with the Hadamard orders for $38\le N<104$, this gives a signing satisfying \eqref{intro:threshold} for every $N\ge38$ and all $p_1,p_2\in[2,\infty]$.

The remaining dimensions are eliminated by comparing explicit sign matrices with lower bounds from order five. Interpolation and congruences of sign-vector images reduce the list to maximal dimensions $3,5,6,9$; sharper estimates remove $3,6,9$. Finally, the minimization in each of the surviving sizes uses symmetries of signs, the geometry of even subsets of a five-element set, and convexity of row multiplicities. The successive reductions are
\[
\begin{gathered}
 \text{infinitely many sizes}\ \longrightarrow\ \text{17 maximal dimensions }(N\le37)\\[3pt]
 \longrightarrow\quad N\in\{3,5,6,9\}\quad\longrightarrow\quad N=5\\[3pt]
 \longrightarrow\quad 5\times5,\ 5\times4,\ 5\times3\ \text{and transposes}\\[3pt]
 \longrightarrow\quad \textcolor{black}{\boxed{\text{five sign matrices and their transposes}.}}
\end{gathered}
\]

\begingroup
For the remaining sizes, the spectral norms, the quadratic moments and absolute sums of the sign images give the required upper bounds by convexity in the reciprocal exponents. In the three sizes with maximal dimension five, the classification gives six representatives. A comparison of the norms of two of the rectangular representatives then reduces the list to five.
\endgroup
Throughout, all sequence spaces and matrices are real unless otherwise indicated. We write $\mathbb N=\{1,2,\ldots\}$ and $t'$ for the conjugate exponent of $t\in[1,\infty]$. For $N\in\mathbb N$, let $B_t^N$ denote the closed unit ball of $\ell_t^N$. Given $A\in\mathbb R^{m\times n}$ and $1\le a,b\le\infty$, we regard $A$ as an operator from $\ell_a^n$ to $\ell_b^m$ and write
\[
 \norm{A}_{a\to b}=\sup_{x\in B_a^n}\norm{Ax}_b.
\]
By duality this is the norm of the bilinear form $(x,y)\mapsto y^tAx$ on $\ell_a^n\times\ell_{b'}^m$, where $A^t$ denotes the transpose of $A$.
\begingroup
Two sign matrices of the same size are called equivalent if one is obtained from the other by permutations and sign changes of rows and columns. Equivalent matrices have the same $\ell_a\to\ell_b$ norm for all $a,b\in[1,\infty]$.
\endgroup

\section{Main results}\label{sec:continuous}
Write $C(p,q)=C_{p,q}$. For a sign matrix $A\in\{-1,1\}^{m\times n}$, define
\begin{equation}\label{eq:normalized}
 F_A(p,q)=\frac{\norm{A}_{p\to q'}}{n^{1/2-1/p}m^{1/2-1/q}\sqrt{\max(m,n)}}.
\end{equation}
Thus $C(p,q)=\sup_{m,n\ge1}\min_A F_A(p,q)$, and $C(p,q)=C(q,p)$ by transposition.

\begingroup
In the six matrices below, $+$ and $-$ denote $1$ and $-1$, respectively. The first classification produces all six candidates; the comparison in Proposition~\ref{prop:D3} removes $D_3$ from the norm minimization.
\endgroup

\[
E_1=\begin{pmatrix}+&+&+&+&+\\+&+&-&+&-\\+&-&-&-&+\\+&+&-&-&+\\+&-&+&+&+\end{pmatrix},\qquad E_2=\begin{pmatrix}+&+&+&+&+\\+&+&+&-&-\\+&+&-&+&-\\+&+&-&-&+\\+&-&+&+&+\end{pmatrix}
\]
\[
D_1=\begin{pmatrix}+&+&+&+\\+&+&-&-\\+&-&+&-\\+&-&-&+\\+&+&+&+\end{pmatrix},\quad D_2=\begin{pmatrix}+&+&+&+\\+&+&+&-\\+&-&-&+\\+&+&-&+\\+&-&+&+\end{pmatrix},\quad D_3=\begin{pmatrix}+&+&+&+\\+&-&-&-\\+&+&+&-\\+&+&-&+\\+&-&+&+\end{pmatrix}
\]
\[
 B=\begin{pmatrix}+&+&+\\+&+&-\\+&-&+\\+&-&-\\+&+&+\end{pmatrix}.
\]
\begingroup
Put $\mathcal E=\{E_1,E_2\}$ and $\mathcal D=\{D_1,D_2\}$.
\endgroup

\begin{theorem}[\textcolor{black}{Reduction to five sign matrices}]\label{thm:continuous}

For every $p,q\in[2,\infty]$,
\begin{equation}\label{eq:band}
\boxed{\begin{aligned}
 C(p,q)=\max\big\{&\min_{A\in\mathcal E}F_A(p,q),\min_{A\in\mathcal D}F_A(p,q),\\
 &\min_{A\in\mathcal D}F_{A^t}(p,q),\ F_B(p,q),\ F_{B^t}(p,q)\big\}.
\end{aligned}}
\end{equation}
The minima over all sign matrices of sizes $5\times5$, $5\times4$ and $5\times3$ are attained, respectively, in $\mathcal E$, in $\mathcal D$ and at $B$. No dimension other than these three and their transposes attains the supremum defining $C(p,q)$.
\end{theorem}

\begingroup
In \eqref{eq:band} all operator norms are taken between spaces of dimension at most five.
\endgroup

For the boundary formulas write $E=E_1$ and, for $1\le r\le2$, set
\begin{equation}\label{eq:LR}
 L(r)=\frac{1}{\sqrt5}\left(\frac{5^r+3^r+3}{5}\right)^{1/r},\qquad
 R(r)=\frac{1}{\sqrt3}\left(\frac{2\cdot3^r+3}{5}\right)^{1/r}.
\end{equation}
\begin{theorem}[Sharp boundary constants]\label{thm:boundary}
For every $p\in[2,\infty]$,
\begin{align}
 C(p,2)=C(2,p)&=\frac{\norm{E}_{p\to2}}{5^{1-1/p}},\label{eq:edge2}\\
 C(\infty,p)=C(p,\infty)&=\max\{L(p'),R(p')\}.\label{eq:edgeinfty}
\end{align}
In particular,
\begin{equation}\label{eq:sharp-corners}
 \boxed{\begin{gathered}
 C(2,2)=\sqrt{\frac85},\qquad C(2,\infty)=C(\infty,2)=\frac{\sqrt{37}}5,\\[3pt]
 C(\infty,\infty)=\frac{3\sqrt3}{5},\qquad
 \sup_{p,q\in[2,\infty]}C(p,q)=\sqrt{\frac85}.
 \end{gathered}}
\end{equation}
\end{theorem}

\begin{proposition}\label{prop:regions}
For the indicated values of $p$ and $q'$,
\begin{equation}\label{eq:regionsB}
\begin{array}{c@{\qquad}c@{\qquad}c}
 p & \text{range of }q' & C(p,q)\\ \midrule
 p\ge24 & [7/4,19/10] & F_E(p,q)\\
 p\ge50 & [3/2,153/100] & F_B(p,q)\\
 p\ge100 & [5/4,3/2] & F_B(p,q)\\
 p\ge400 & [6/5,5/4] & F_B(p,q)
\end{array}
\end{equation}
The transposed identities hold with $p,q$ interchanged and the matrices transposed.
\end{proposition}

\subsection*{Two examples}
\emph{The pair $(\infty,4)$.} Here $r=4/3$. Theorem~\ref{thm:boundary} gives
\[
 L(4/3)^{4/3}=\frac{5\sqrt[3]{5}+3\sqrt[3]{3}+3}{5\sqrt[3]{25}}
 <\frac{11}{10}<\frac{6\sqrt[3]{3}+3}{5\sqrt[3]{9}}=R(4/3)^{4/3}.
\]
Indeed, $\sqrt[3]{5}<1.71$, $1.44<\sqrt[3]{3}<1.45$, $\sqrt[3]{25}>2.9$ and $\sqrt[3]{9}<2.1$, as is seen by cubing. The matrix $B$ therefore realizes the constant in size $5\times3$. At $u=(1,1,1)^t$, of $\ell_\infty^3$-norm one, $Bu=(3,1,1,-1,3)^t$. Its absolute values form a maximizing profile for $B$ (see Section~\ref{sec:classification}). Consequently
\begin{equation}\label{eq:example-boundary}
 \boxed{C(\infty,4)=C(4,\infty)=\frac{1}{\sqrt3}
 \left(\frac{2\cdot3^{4/3}+3}{5}\right)^{3/4}\approx1.08907.}
\end{equation}
For the transposed pair $(4,\infty)$ the realizing matrix is $B^t$.

\pagebreak[2]
\noindent\emph{The pair $(100,5)$.} Since $q'=5/4$, Proposition~\ref{prop:regions} gives
\begin{equation}\label{eq:example-interior}
 C(100,5)=\frac{\|B\|_{100\to5/4}}{3^{49/100}5^{4/5}}.
\end{equation}
Testing on the normalized constant vector gives $F_B(100,5)\ge R(5/4)$. On the other hand, $B^t B$ has largest eigenvalue $7$, since $B^tB=4I_3+J_3$, where $J_3$ is the all-ones matrix, so $F_B(2,5)\le F_B(2,2)=\sqrt{7/5}$. Interpolation in $1/p$ yields
\[
 \underbrace{R(5/4)}_{\approx1.0766913}\le C(100,5)\le
 \underbrace{(7/5)^{1/100}R(5/4)^{49/50}}_{\approx1.0787248}.
\]

\section{Normalized norms and the small sign matrices}\label{sec:classification}
Write $p_t=2/(1-t)$ for $0\le t<1$, with $p_1=\infty$. For $u\in\mathbb R^k$, let $M_s(u)=k^{-1/s}\|u\|_s$ for $1\le s<\infty$ and $M_\infty(u)=\|u\|_\infty$. In particular,
\[
 F_A(p,q)=\frac{1}{\sqrt{mnN}}\sup_{M_p(u)\le1,\,M_q(v)\le1}|v^tAu|,\qquad N=\max(m,n).
\]
The monotonicity of the power means proves that $F_A$ is nonincreasing in each exponent. Also, if $x=(1-\alpha)x_0+\alpha x_1$ with $0\le\alpha\le1$, then $1/p_x'=(1-\alpha)/p_{x_0}'+\alpha/p_{x_1}'$. H\"older's inequality gives
\[
 \|A^tv\|_{p_x'}\le \|A^tv\|_{p_{x_0}'}^{1-\alpha}\|A^tv\|_{p_{x_1}'}^\alpha.
\]
Taking the supremum over $v\in B_{p_y}^m$ and dividing by the dimensional factors in \eqref{eq:normalized}, which are exponential affine functions of $x$, proves convexity of $\log F_A(p_x,p_y)^2$ in $x$. Duality gives the analogous convexity in $y$. If $a,b,c,d$ bound the squares of $F_A$ at $(2,2),(\infty,2),(2,\infty),(\infty,\infty)$, respectively, and $t=\min(x,y)$, these facts give
\begin{equation}\label{eq:corners}
 F_A(p_x,p_y)^2\le a^{(1-t)^2}(bc)^{t(1-t)}d^{t^2}
 \le a(1-t)^2+(b+c)t(1-t)+dt^2.
\end{equation}
For a sign vector $u\in\{-1,1\}^n$, put $z=|Au|$ coordinatewise. Direct substitution of $n^{-1/p}u$ gives
\begin{equation}\label{eq:profile}
 F_A(p,q)^2\ge \frac{m}{nN}\left(\frac1m\sum_i z_i^{q'}\right)^{2/q'},\qquad N=\max(m,n).
\end{equation}
The transposed estimate interchanges $p,q$. A profile will mean the absolute values of the image of a sign vector, with zeros included. At $p=\infty$ the norm is the maximum over profiles, since a convex function on the cube attains its maximum at a vertex.

\begin{lemma}\label{lem:tangent}
For a nonzero nonnegative vector $z$ of length $M$, with $S=\sum z_i$, and $\gamma>0$,
\[
 \gamma\left(\frac1M\sum_i z_i^{2/(1+t)}\right)^{1+t}
 \ge \gamma(S/M)^2 e^{k(1-t)},\qquad
 k=\frac{\sum_i z_i\log z_i}{S}-\log(S/M).
\]
\end{lemma}
\begin{proof}The logarithm of the left side is convex, being the perspective of log-sum-exp on the positive coordinates. Its derivative at $t=1$ is $-k$. Use the tangent there, with $0\log0=0$.\end{proof}

\begingroup
\begin{lemma}[Balancing multiplicities]\label{lem:balance}
Let $v_1,\ldots,v_k\in\{-1,1\}^n$ represent row types modulo sign. Suppose a signed permutation of the columns interchanges types $v_i,v_j$ and preserves the other multiplicities. For $c\in\mathbb R_+^k$ write
\[
 \Psi(c)=\sup_{x\in B_p^n}\sum_{\nu=1}^k c_\nu|\langle v_\nu,x\rangle|^{q'}.
\]
If $c_i,c_j$ are integers with $c_i\ge c_j+2$, replacing $(c_i,c_j)$ by $(c_i-1,c_j+1)$ does not increase $\Psi$.
\end{lemma}
\begin{proof}
The function $\Psi$ is convex as the supremum of linear functions of $c$. The assumed symmetry gives $\Psi(c)=\Psi(c^*)$, where $c^*$ is obtained by interchanging the two multiplicities. With $\alpha=1/(c_i-c_j)$, the new multiplicity vector is $(1-\alpha)c+\alpha c^*$. Convexity gives the assertion. If all transpositions of the types are induced by column symmetries, iteration terminates at multiplicities differing by at most one, since the sum of their squares decreases by $2(c_i-c_j-1)>0$ at each transfer.
\end{proof}
For integer multiplicities, $\Psi(c)$ is the $q'$-th power of the norm of the matrix formed by repeating the row $v_i$ exactly $c_i$ times.
\endgroup
For three columns, normalize each row by a sign so that its product is one. The four types are the constant row and the three rows with two negative entries. Coordinate permutations fix the constant row and induce all permutations of the other three. Column sign changes of product one act transitively on the four types. Thus the group induced on the types is $S_4$. \begingroup
Lemma~\ref{lem:balance} gives the multiplicities $(2,1,1,1)$ when there are five rows. Hence $B$ minimizes the norm in size $5\times3$.
\endgroup Its profile $(3,3,1,1,1)$ is attained by a sign vector equal to the repeated row. More generally every $5\times3$ matrix has such a lower bound since two rows coincide up to sign. Thus $C(p,q)\ge R(q')$ and $C(p,q)\ge R(p')$.

\subsection{The square matrices of order five}
Normalize the product of the five entries of each row to one, and identify the row with its even subset of negative positions in $\{1,2,3,4,5\}$. Even column sign changes act by symmetric difference, and column permutations act by permutations of the ground set. Join distinct rows when their subsets differ in four positions, so their inner product is $-3$; all other distinct rows have inner product $1$.
\begin{lemma}[\textcolor{black}{Classification by even subsets}]\label{lem:graph}
Every $5\times5$ sign matrix is equivalent to $E_1$ or $E_2$, or has a profile whose $r$-power sum is at least $5^r+2\cdot3^r+2$ for all $1\le r\le2$. The classes of $E_1,E_2$ are invariant under transposition.
\end{lemma}
\begin{proof}
If two rows coincide, evaluate on one of them, giving at least $(5,5,1,1,1)$, and use $5^r+1\ge2\cdot3^r$. A vertex of degree at least two yields a profile dominating $(5,3,3,1,1)$. Otherwise the graph of five distinct rows is a matching.

With no edges, translate a row to the empty set. The other four are distinct intersecting pairs and have a common element (a triangle of pairs cannot be extended to a fourth pair meeting all three). The four-set omitting this element has inner product $-3$ with all five rows. Its profile $(3,3,3,3,3)$ dominates $(5,3,3,1,1)$ in $r$-power sum by concavity of $s^{r/2}$ and equality of squared sums.

With one edge, translate an endpoint to $\varnothing$ and the other to $1234$, abbreviating subsets by their elements. The three remaining pairs lie in $1234$ and intersect pairwise. They form a star or triangle, equivalent by translation by $1234$. The representative is $\{\varnothing,1234,12,13,14\}$, and the set $2345$ yields $(3,3,3,3,1)$. With two edges, translate the isolated vertex to $\varnothing$. The four remaining pairs form two disjointness edges and otherwise intersect. Up to permutation they are $12,34,13,24$. This is the class of $E_1$; the one-edge case is $E_2$. Direct normalization of the rows of $E_1,E_2,E_1^t,E_2^t$ gives respectively
\[
(\varnothing,35,15,34,1345),\quad(\varnothing,45,35,34,1345),
\]
\[
(\varnothing,35,15,34,1345),\quad(\varnothing,1234,34,24,23).
\]
Their graphs give the asserted equivalences.
\end{proof}
The column Gram forms of $E_1,E_2$, up to sign changes and permutations, are respectively
\begin{align*}
 Q_2(x)&=4(x_1-x_2)^2+4(x_3-x_4)^2+4x_5^2+(\sum x_i)^2,\\
 Q_1(x)&=4(x_1-x_2)^2+4(x_3^2+x_4^2+x_5^2)+(\sum x_i)^2.
\end{align*}
Both have maximum $37$ on sign vectors. For $Q_2$, the two pair-difference directions have eigenvalue $8$. On their orthogonal complement, $Q_2(u,u,v,v,w)=4w^2+(2u+2v+w)^2$. Since
\[
 (2u+2v+w)^2\le (16u^2+16v^2+4w^2)(4/16+4/16+1/4),
\]
we have $Q_2(u,u,v,v,w)\le16u^2+16v^2+8w^2$. Hence $\|E_1\|_{2\to2}^2=8$.

In the canonical two-edge configuration, an exterior pair is disjoint from at most two of the four pairs, and a four-set is adjacent to the empty set and at most two pairs. Thus a profile with no $5$ has at most three $3$'s, and a profile with a $5$ has at most one $3$. The profile $(5,3,1,1,1)$ is attained and dominates all the others in $r$-power sum by $5^r+1\ge2\cdot3^r$. In the one-edge class the profile $(3,3,3,3,1)$ dominates this profile by concavity of the powers of their squares.

\begin{proposition}\label{prop:square}
For all $p,q\ge2$, the minimum over the $5\times5$ sign matrices is $\min(F_{E_1},F_{E_2})$.
\end{proposition}
\begin{proof}The four corner squares for $E_1$ are $8/5,37/25,37/25,121/125$. Thus \eqref{eq:corners} bounds its square by $W=(8/5)(1-t)^2+(74/25)t(1-t)+(121/125)t^2$. For any other class, Lemma~\ref{lem:graph} applies in either orientation. Lemma~\ref{lem:tangent} with $(5,3,3,1,1)$ gives a lower bound
\[
 H=\frac{169}{125}\left(1+\frac{17}{100}(1-t)+\frac{289}{20000}(1-t)^2\right),
\]
since $k=(5\log5+6\log3)/13-\log(13/5)>17/100$. In fact $H-W=3(342947t^2-24094t+1147)/2500000>0$ by the negative discriminant.\end{proof}

\begingroup
\begin{corollary}\label{cor:Eregion}
If $p,q\ge4$, then $F_{E_1}(p,q)<F_{E_2}(p,q)$.
\end{corollary}
\begin{proof}
Put $t=\min(1-2/p,1-2/q)$, so $t\ge1/2$. From the four corners for $E_1$ and \eqref{eq:corners}, using $\log(8/5)<.471$, $\log(37/25)<.393$ and $\log(121/125)<-.032$, we get
\[
 \log F_{E_1}(p,q)^2\le .471-.156t-.347t^2.
\]
The profile $(3,3,3,3,1)$ of $E_2$, in either orientation, gives by the tangent at $t=0$ a lower bound $.3920-.1364t$ (Section~\ref{sec:arithmetic}). The difference between this lower bound and the upper bound for $E_1$ is
\[
 -.079+.0196t+.347t^2\ge -.079+.0098+.08675=.01755>0.
\]
This proves the strict comparison.
\end{proof}
\endgroup
\subsection{The five-by-four matrices}
There are eight row types up to sign, divided into two parity families of four mutually orthogonal types, with cross inner products of modulus two. The same convexity of multiplicities used for three columns applies. Within either parity family the induced permutation group is $S_4$: in the odd family take representatives with exactly one negative coordinate, and in the even family conjugate these symmetries by a column sign change of odd parity. For a parity split $5+0$, balance the four types to $(2,1,1,1)$, giving $D_1$. For $4+1$, normalize the minority row to the constant row; permutations of coordinates induce all permutations of the four odd types, represented by one negative coordinate. Balancing gives $D_3$.

For $3+2$, the Klein four-group of double transpositions fixes each even type up to sign and is transitive on the odd types. \begingroup
If the minority rows coincide, choose a double transposition sending their type to a different odd type. It preserves the multiplicities of the even types. Averaging the two multiplicity vectors with minority multiplicities $(2,0)$ and $(0,2)$ gives $(1,1)$; convexity of $\Psi$ proves that the norm cannot increase.
\endgroup If the majority types are distinct, normalize the missing even type to the constant row. All coordinate permutations preserve the three majority types as a set and move the minority pair to any other pair; this is $D_2$.

In the remaining case normalize the repeated majority row to the constant type. Evaluation on that row yields a profile dominating $(4,4,2,2,0)$. In the transpose, choose minority rows with single negatives in distinct positions. For three identical majority rows the all-ones test gives $(5,5,3,3)$; for precisely two, choose the sign of the third majority row so its two negative positions meet the minority negative positions. The profile is then $(5,5,1,1)$ or $(5,3,3,1)$, the former dominating the latter in $r$-power sum. The corner squares for $D_1$ are $8/5,8/5,28/25,1$: this follows from the four orthogonal Hadamard rows and a repetition. In the transposed orientation the Walsh transform of four signs has one entry of modulus four or four of modulus two, and the repeated row adds a sign to every coordinate. Consequently
\[
 F_{D_1}^2\le 8/5-(12/25)t-(3/25)t^2.
\]
Both discarded profiles have normalized square $36/25$ at $t=1$, with $k>1/9$ in Lemma~\ref{lem:tangent} (the values are $(5/3)\log2-\log(12/5)$ and $(5\log5-6\log3)/12$). Their lower bound $\frac{36}{25}(1+(1-t)/9)=8/5-4t/25$ exceeds the upper bound for $D_1$, strictly also at $t=0$ because $k>1/9$. This proves the minimum over $5\times4$ is attained among $D_1,D_2,D_3$.

\begingroup
\begin{proposition}\label{prop:D3}
For all $p,q\in[2,\infty]$, one has
\[
 F_{D_1}(p,q)\le F_{D_3}(p,q),\qquad
 F_{D_1^t}(p,q)\le F_{D_3^t}(p,q).
\]
Both inequalities are strict when the exponent corresponding to the five-dimensional factor is greater than two. Consequently the minimum in size $5\times4$ is attained at $D_1$ or $D_2$.
\end{proposition}
\begin{proof}
Set $y=1-2/q$. By monotonicity in $p$, and convexity of the logarithm of the squared normalized norm in $y$,
\[
 F_{D_1}(p,q)^2\le F_{D_1}(2,q)^2
 \le (8/5)^{1-y}(28/25)^y=\frac85(7/10)^y.
\]
On the other hand the sign vector $(1,1,1,1)^t$ gives the absolute profile $(4,2,2,2,2)$ for $D_3$. With $r=2/(1+y)$, \eqref{eq:profile} gives
\[
 F_{D_3}(p,q)^2\ge M(y):=\frac14\left(\frac{4^r+4\cdot2^r}{5}\right)^{1+y}.
\]
The function $\log M$ is convex by the perspective of log-sum-exp. At zero,
\[
 M(0)=8/5,\qquad (\log M)'(0)=\log(32/5)-2\frac{16\log4+16\log2}{32}=\log(4/5).
\]
Its supporting line at zero therefore gives $M(y)\ge(8/5)(4/5)^y$. Since $7/10<4/5$, the claimed inequality follows, strictly for $y>0$. The transposed comparison follows by interchanging $p,q$. The fixed-size classification just proved now gives the last assertion.
\end{proof}
\endgroup
\begingroup
Both $D_1$ and $D_2$ are needed in the minimum for the fixed size $5\times4$. Indeed, the four corner bounds for $D_1$ are equalities, so $F_{D_1}(\infty,2)^2=8/5$ and $F_{D_1}(2,\infty)^2=28/25$. For $D_2$, the two parity families of three and two distinct types yield a largest direct sign-profile energy of $28$ (the profile $(4,2,2,2,0)$), so $F_{D_2}(\infty,2)^2=7/5<8/5$. Also $D_2^t\mathbf1=(5,1,1,3)^t$, so $F_{D_2}(2,\infty)^2\ge36/25>28/25$. These fixed-size comparisons do not prove that either rectangular matrix is indispensable for the global maximum in \eqref{eq:band}.
\endgroup
Let $C_5(p,q)$ denote the maximum of the normalized minima over $5\times5$, $5\times3$ and $3\times5$. The row-profile bounds for $E_1,E_2$ and $B$, together with \eqref{eq:profile}, give, for $r=2/(1+t)$,
\begin{equation}\label{eq:lower}
 C(p,q)\ge C_5(p,q)\ge\max\{L(r),R(r)\},\qquad C_5(p,q)^2\ge27/25.
\end{equation}
For the first bound use the orientation corresponding to the smaller exponent and \eqref{eq:profile}; the square classes are invariant under transposition.

\section{Hadamard orders and the first reduction}\label{sec:had}
A submatrix of an order-$h$ Hadamard matrix satisfies $F_A^2\le h/N$, by extension with zeros and the inclusions of the normalized balls. Over a finite field of odd order $P$, the quadratic character $\chi$ has $\chi(0)=0$, $\sum_x\chi(x)=0$, and $\sum_x\chi(x-a)\chi(x-b)=-1$ if $a\ne b$, by the substitution $(x-a)/(x-b)$. Set $T_{ij}=\chi(j-i)$. If $P\equiv3\pmod4$, $T$ is skew symmetric, $TT^t=PI-J$, and $\left(\begin{smallmatrix}1&\mathbf1^t\\-\mathbf1&I+T\end{smallmatrix}\right)$ is Hadamard of order $P+1$. If $P\equiv1\pmod4$, the symmetric conference matrix $K=\left(\begin{smallmatrix}0&\mathbf1^t\\\mathbf1&T\end{smallmatrix}\right)$ satisfies $K^2=PI$, and $\left(\begin{smallmatrix}K+I&K-I\\K-I&-K-I\end{smallmatrix}\right)$ is Hadamard of order $2(P+1)$. Tensor products preserve the Hadamard property. This supplies the orders
\[
1,2,4,8,12,16,20,24,28,32,36,40,44,48,52,56,60,64,
\]
\[68,72,76,80,84,88,96,100,104,180.\]
The Paley orders of the first type use $P=11,19,23,43,59,67,83,179$; the second type uses $P=13,17,25,37,49$. The other orders are powers of two or doublings. For $N\ne5$, $N<38$, these give $h/N\le4/3$. For $38\le N\le103$, round up to the next multiple of four except for $89\le N\le92$, when $h=96$. The largest ratio is at most $96/89<27/25$. For $N\ge104$, take $j\ge0$ with $104\le N/2^j<208$ and round up to an order from
\[
104,112,120,128,136,144,152,160,168,176,180,192,200,208.
\]
Consecutive ratios are at most $14/13$, so tensoring by the order-$2^j$ Sylvester matrix gives $h/N\le14/13<27/25$. By \eqref{eq:lower}, all $N\ge38$ are eliminated. For $3\le N\le37$ only the orders with $4\lceil N/4\rceil/N\ge27/25$ need be retained; orders one and two have bound one. Consequently,
\begin{equation}\label{eq:finite-reduction}
 C(p,q)=\max_{\substack{1\le m,n\le37\\ \max(m,n)\in\mathcal N}}
 \min_{A\in\{-1,1\}^{m\times n}}F_A(p,q),
\end{equation}
where $\mathcal N=\{3,5,6,7,9,10,11,13,14,17,18,21,22,25,29,33,37\}$.
\section{Elimination of the larger orders}\label{sec:large}
Set $t=\min(1-2/p,1-2/q)$ and $r=2/(1+t)$. For a witness of size $N\times n$ obtained from a Hadamard matrix of order $h$, let $E,G,L_0$ bound $\|A\|_{\infty\to2}^2,\|A^t\|_{\infty\to2}^2,\|A\|_{\infty\to1}$, respectively. Equation \eqref{eq:corners} gives
\begin{equation}\label{eq:W}
 F_A^2\le W(t)=a(1-t)^2+bt(1-t)+dt^2,\quad
 a=h/N,\quad b=E/(Nn)+G/N^2,\quad d=L_0^2/(N^2n).
\end{equation}
Put $f(t)=R(2/(1+t))^2$. Its logarithm is convex, so $f$ is convex. \begingroup
At $t=1$, $f(1)=27/25$ and $f'(1)=-(9/25)\log(125/81)<-.156$. At $t=1/2$, $1.1860<f(1/2)<1.1862$ and $-.287<f'(1/2)<-.286$.
\endgroup
 The tangents give
\begin{equation}\label{eq:H}
 f(t)\ge \textcolor{black}{H_1(t)=\frac{1236-156t}{1000}},\qquad
 f(t)>\textcolor{black}{H_3(t)=\frac{1329-287t}{1000}}.
\end{equation}
The rational bounds follow from the estimates in Section~\ref{sec:arithmetic}.

\subsection{Parity and deletion}
For an $m\times n$ submatrix of an order-$h$ Hadamard matrix and a sign vector $u$, put $z=|Au|$. Then
\begin{equation}\label{eq:energy}
 \sum_i z_i^2\le B_0=\min\{hn-(h-m)(n\bmod2),mn^2\}.
\end{equation}
Round $B_0$ down to $B$ in the class $m\pmod8$ if $n$ is odd, or $0\pmod4$ otherwise. For consecutive nonnegative integers $u_0,v_0=u_0+2$ of parity $n$ with $mu_0^2\le B\le mv_0^2$, $(z_i-u_0)(z_i-v_0)\ge0$ yields
\begin{equation}\label{eq:parity}
 \sum z_i\le (B+mu_0v_0)/(u_0+v_0).
\end{equation}
Round down to parity $mn$; when $B=mn^2$ use $mn$. Denote this bound by $U_0(m,n,h)$. The initial bounds in \eqref{eq:W} are $E=B(N,n,h),G=B(n,N,h),L_0=\min(U_0(N,n,h),U_0(n,N,h))$.

When $n$ is even, let $k,K$ count the positive products of the selected $n$ entries in the retained rows and in all $h$ rows, respectively. For any input sign vector, the numbers of retained coordinates congruent to $0$ and $2$ modulo four, and the number in the latter class among all rows, form one of the triples
\[
(k_0,k_2,K_2)=(k,m-k,h-K),\quad(m-k,k,K).
\]
Indeed the residue is determined by the parity of the number of disagreements. Let $B_*$ be the largest integer congruent to $4k_2\pmod{16}$ not exceeding $\min(hn-4(K_2-k_2),mn^2)$. The omitted coordinates in class two contribute at least four to the energy, and the squares have the asserted congruence. Hence for $s\ge1$,
\begin{equation}\label{eq:mod16}
 \sum z_i^2\le B_*,\qquad \sum z_i\le\frac{B_*+ms^2-k_0\operatorname{dist}(s,4\mathbb Z)^2-k_2\operatorname{dist}(s,4\mathbb Z+2)^2}{2s}.
\end{equation}
The latter follows by expanding $(z_i-s)^2$. Take the minimum over $s=1,\ldots,7$ and \eqref{eq:parity}, round down to $2k_2\pmod4$, and maximize over the two triples; call the result $U(m,n,h;k,K)$. For energy take the maximum of the two $B_*$ and the minimum with $B$.

If $L(A)=\|A\|_{\infty\to1}$, appending a column to $A$ gives $L(A)\le L(A^+)-(m\bmod2)$: extend an optimizing pair of signs, choosing the new sign to make the contribution nonnegative and, for odd $m$, at least one. If $A_{-j}$ deletes column $j$, then $(n-1)L(A)\le\sum_j L(A_{-j})$, by summing the restricted evaluations of an optimizing pair. Both statements hold for transposes.

\subsection{The Walsh construction at order thirteen}
Use $H_{16}=H_4\otimes H_4$ with lexicographic indices from zero and retain the first thirteen rows. The three omitted rows have distinct second coordinates. With twelve columns, take three full blocks of four. The Walsh transform of four input signs has either one coordinate of modulus four (concentrated type), or four coordinates of modulus two (flat type), whose signs have product negative. This follows by distinguishing the parity of the input; the even sign vectors are rows of $H_4$ up to sign, and odd vectors are equivalent to a vector with one negative. At each second-coordinate frequency, apply the outer Walsh transform to three coefficients and a zero. Its absolute sum is $4\max\{\max(|a|,|b|,|c|), (|a|+|b|+|c|)/2\}$, as follows by summing the four sign combinations modulo global sign. For thirteen columns add the first column of the fourth block. Its four frequency signs have positive product. The following bounds result (the number of flat blocks is denoted by $f_0$).
\begin{center}\small
\begin{tabular}{clrrrr}\toprule
 & &\multicolumn{2}{c}{$n=12$}&\multicolumn{2}{c}{$n=13$}\\
$f_0$& Concentrated frequencies & full & retained & full & retained\\\midrule
3& none &48&42&52&45\\
2& one &40&40&48&45\\
1& same &44&38&46&43\\
1& distinct &48&42&48&45\\
0& three distinct &48&40&52&45\\
0& at most two &32&32&42&39\\\bottomrule
\end{tabular}
\end{center}
\begingroup
The absolute-sum identity for three coefficients follows from sorting their absolute values as $a\ge b\ge c\ge0$: the four absolute values are $a+b+c$, $a+b-c$, $a-b+c$, and $|a-b-c|$. Their sum is $4a$ if $a\ge b+c$ and $2(a+b+c)$ otherwise. For twelve columns the frequency triples in the six cases in the table are respectively four copies of $(2,2,2)$; $(4,2,2)$ and three copies of $(0,2,2)$; $(4,4,2)$ and three copies of $(0,0,2)$; two copies of $(4,0,2)$ and two of $(0,0,2)$; three copies of $(4,0,0)$ and one zero triple; or at most two nonzero triples of total absolute sum twelve. Substitution in the identity gives the full bounds $48,40,44,48,48,32$ in the displayed order.
\endgroup
For $n=12$, one or three flat blocks imply all coordinates have modulus at least two. Three distinct concentrated frequencies yield three groups of modulus four and at least two of the omitted coordinates belong to these. For $n=13$ the four-coefficient absolute sums are bounded as follows:
\[
\begin{array}{c|rrrrrrrrr}
\text{moduli}&2221&4221&2201&4421&4201&2001&4001&4401&4441\\\hline
\text{sum}&14&18&10&22&16&8&16&18&26
\end{array}
\]
These values follow from the even and odd sign combinations; a zero makes the parity immaterial. For two flat blocks the full sum is at most $18+3\cdot10=48$. For one flat block it is at most $22+3\cdot8=46$ when the concentrated frequencies coincide, and $2\cdot16+2\cdot8=48$ otherwise. With no flat blocks the corresponding bounds are $3\cdot16+4=52$ for three frequencies, $18+16+2\cdot4=42$ for two frequencies and $26+3\cdot4=38$ for one frequency. In particular $4421$ yields profiles $(11,5,1,1)$ or $(9,7,3,3)$, giving $18$ or $22$. For $2221$ the profiles are $(7,1,1,1)$ or $(5,3,3,3)$; the product of the four frequency sign products is negative, so there are one or three groups of the latter type. In the first case the total is $44$; in the second case at least two omitted coordinates have modulus at least three, so at least seven is lost. With three concentrated distinct frequencies the three active groups have modulus at least three, again losing at least seven. In all other cases the odd coordinates lose at least three. This proves $L_{13,12}\le42$, $L_{13,13}\le45$ and, by deletion, $L_{13,11}\le41$.

\begingroup
\begin{lemma}\label{lem:paley51}
Let $A$ be the leading $17\times9$ submatrix of the order-twenty Paley matrix (first construction in Section~\ref{sec:had}, with field elements in the order $0,1,\ldots,18$). Then $\|A\|_{\infty\to1}=51$.
\end{lemma}
\begin{proof}
Normalize each row of $A$ to have product one, obtaining $A_0$. This leaves the norm unchanged. A sign vector $u\in\{-1,1\}^9$ may also be assumed to have product one, by replacing it by $-u$ if needed. All signed coordinates of $A_0u$ belong to $\{9,5,1,-3,-7\}$. Their squared sum is at most $177$: the full Hadamard image has energy $20\cdot9=180$, and the three omitted odd coordinates contribute at least three. The absolute sum is odd, and the inequality $(z-3)(z-5)\ge0$ for positive odd $z$ gives
\[
 8\sum_{i=1}^{17}|(A_0u)_i|\le177+15\cdot17=432.
\]
Thus the absolute sum is at most $53$.

Suppose it equals $53$. Denote by $b,c,d,e$ the multiplicities of the absolute values $3,5,7,9$. The first and second moments give
\[
 b+2c+3d+4e=18,\qquad b+3c+6d+10e\le20.
\]
Hence $c+3d+6e\le2$, so $d=e=0$. The length condition $b+c\le17$ implies $c\ge1$. There are exactly two possibilities: $(b,c)=(16,1)$ or $(14,2)$, with respectively zero or one remaining coordinate of modulus one. In the signed normalization their sums are $-43$ and $-31$, respectively.

For the Paley construction the nonzero squares modulo $19$ are $1,4,5,6,7,9,11,16,17$. Substituting these residues in the entries and multiplying each row by its product gives the column sums
\[
 s=\mathbf1^t A_0=(-1,5,-3,3,-7,1,-5,3,5).
\]
In particular $\|s\|_1=33$, so $s u=-43$ is impossible. If $su=-31$, the vector $u$ differs from $-\operatorname{sign}(s)$ in precisely one coordinate of weight one: either the first or the sixth. For these two vectors, respectively,
\begin{align*}
 u&=(-1,-1,1,-1,1,-1,1,-1,-1)^t,\\
 A_0u&=(-3,1,1,-3,1,-3,-3,1,1,1,1,-3,-7,-3,-7,-3,-3)^t;\\
 u&=(1,-1,1,-1,1,1,1,-1,-1)^t,\\
 A_0u&=(1,1,-3,-7,1,-3,1,1,1,-3,-3,-7,-3,-3,-3,-3,1)^t.
\end{align*}
Each has a coordinate of modulus seven, contradicting $d=0$. Thus the norm is at most $51$.
Finally, the vector $(-1,1,-1,-1,-1,1,-1,1,1)^t$ has image
\[
 (-1,-5,-1,3,3,3,3,-5,3,3,3,3,3,3,3,3,3)^t,
\]
whose absolute sum is $51$. This proves equality.
\end{proof}
\endgroup
\begingroup
\begin{lemma}\label{lem:paley117}
Let $A$ be the leading $25\times21$ submatrix of the order-twenty-eight Hadamard matrix in the second Paley construction of Section~\ref{sec:had}, with field elements in the order $0,1,\ldots,12$. Then $\|A\|_{\infty\to1}\le117$.
\end{lemma}
\begin{proof}
Normalize the rows to product one, obtaining $A_0$, and assume that $u\in\{-1,1\}^{21}$ has product one. The signed coordinates of $A_0u$ are congruent to $1$ modulo four, their absolute values $z_i$ are odd, and their energy is at most $585=28\cdot21-3$. The inequality $(z_i-3)(z_i-5)\ge0$ gives $8\sum z_i\le585+375=960$. The sum is odd and hence at most $119$.

If the sum is $119$, then $\sum (z_i-3)(z_i-5)\le8$. Each nonzero summand is at least eight, with equality only at $1$ and $7$. There are therefore just three possible profiles: $(5^{22},3^3)$, $(5^{23},3,1)$, or $(5^{20},3^4,7)$, with powers denoting multiplicities. In the signed normalization the coordinate sums are respectively $101,113,81$.

The nonzero squares modulo $13$ are $1,3,4,9,10,12$. The conference construction and row normalization give the column sums
\[
 s=(3,-3,1,-11,-9,5,-3,3,1,5,-3,3,5,3,-3,3,-1,-5,-7,-5,3),
\]
with $\|s\|_1=85$. Thus only the sum $81$ is possible. The vector $u$ must differ from $\operatorname{sign}(s)$ in exactly two of the three coordinates of weight one (indices $2,8,16$, starting from zero). Indeed, the sum of the weights of the changed coordinates is $(85-81)/2=2$, and all other weights are odd and at least three. The normalized row with index $22$ is
\[
 (1,-1,-1,-1,1,1,-1,1,-1,1,-1,1,1,-1,-1,1,1,1,-1,-1,1).
\]
Its inner product with $\operatorname{sign}(s)$ is $9$, and its contributions at the three specified positions are all $-1$. Changing any two of them makes the inner product $13$. This is impossible in the three profiles above. Thus $119$ cannot occur and the absolute sum is at most $117$.
\end{proof}
\endgroup
\subsection{The other witnesses and scalar comparisons}
Use leading submatrices of Sylvester matrices for orders $8,16,32$, Paley matrices of the first type for $12,20,24$, of the second type for $28,36$, and $H_2\otimes H_{20}$ for $40$. In a Paley matrix of order $h=P+1$ the first row is positive, the first column is negative outside its first entry, and the remaining entry is positive on the diagonal and for $j-i$ a nonzero quadratic residue. At order 28 use the symmetric conference construction with $P=13$. The following exceptional choices use zero-based indices: for $13\times9$ choose columns $0,1,2,3,4,5,8,10,12$; for $17\times17$, columns $0,\ldots,15,19$; for $21\times21$, use $H_2\otimes H_{12}$ and omit rows and columns $0,1,12$.

The bounds for $L(A)$ obtained from \eqref{eq:mod16} are listed below. A T indicates transposition before applying \eqref{eq:mod16}; $+j$ appends column $j$ before applying that inequality and then uses the extension bound. The counts $k,K$ are the row-product counts from the given matrices, and the final column is the resulting bound for $L(A)$.
\begin{center}\small
\begin{tabular}{rrcrrr@{\qquad}rrcrrr}\toprule
$N$&$n$&op.&$k$&$K$&$L_0$&$N$&$n$&op.&$k$&$K$&$L_0$\\\midrule
10&7&T&3&6&26&17&10&--&12&14&54\\
10&9&T&4&6&30&17&11&+16&7&8&57\\
13&4&--&13&16&26&17&12&--&6&8&60\\
13&7&+7&13&16&33&17&14&--&12&14&66\\
13&8&--&13&16&34&17&15&+15&11&12&67\\
13&9&+15&13&16&39&17&16&--&11&12&68\\
13&10&--&7&8&42&21&4&--&11&12&34\\
14&9&T&5&8&42&21&12&--&16&16&74\\
14&14&--&7&8&50&21&14&--&11&12&78\\
17&4&--&11&12&30&21&16&--&15&16&84\\
17&7&+7&11&12&45&21&20&--&15&16&96\\
25&4&--&8&8&42&25&14&--&14&14&92\\
25&16&--&11&12&98&25&23&+23&9&12&123\\
25&25&+25&13&14&129&&&&&&\\\bottomrule
\end{tabular}
\end{center}
For instance, in the $13\times10$ case the product of the chosen Sylvester columns equals the column with binary index $1$, so $(k,K)=(7,8)$ and \eqref{eq:mod16} gives $E\le156$. For $17\times17$, deleting a row and transposing gives the following counts and bounds:
\[
\begin{array}{c|rrrrrrr}
(k,K)&(6,8)&(11,12)&(9,12)&(8,8)&(7,8)&(10,12)&(4,4)\\
\text{multiplicity}&5&4&3&1&2&1&1\\
U&68&68&68&68&70&70&72
\end{array}
\]
Thus $16L\le13\cdot68+3\cdot70+72=1166$, and $L\le71$ by parity. For the tensorial $21\times21$ matrix the same operation gives $(10,12),(11,12),(9,12),(21,24)$ with multiplicities $10,5,5,1$ and $U\le96$ in every case, so $20L\le21\cdot96$, hence $L\le99$. The counts can also be obtained from the fact that the product of columns $0,1,12,13$ in $H_2\otimes H_{12}$ is constant; after the deletion the product is constant or the product of two distinct columns.

\begingroup
The row-product counts can be calculated directly by multiplying the columns specified above. More explicitly, for the index set $J$ used in a witness, define $\eta_i=\prod_{j\in J}H_{ij}$, where $H$ is the displayed Paley or Sylvester construction. Then $k=(|I|+\sum_{i\in I}\eta_i)/2$ and $K=(h+\sum_{i=0}^{h-1}\eta_i)/2$. For a transposed witness interchange $I,J$ and use $H^t$; for an appended column include it in $J$ before calculating the products. In the Sylvester case $\eta_i=(-1)^{i\cdot(\bigoplus_{j\in J}j)}$, where the scalar product and the sum of binary indices are over $\mathbb F_2$. In the first Paley case the product is obtained by counting nonresidues among the differences $j-i$ for the selected columns (with the first column treated as specified above). The second construction uses the four displayed blocks $K+I,K-I,K-I,-K-I$. These formulas give the counts without examining any optimizing sign vector.

For the deletion calculations, the following are the deleted row indices producing each count pair (indices start at zero):
\[
\begin{array}{c|l|l}
N&(k,K)&\text{deleted row indices}\\\hline
17&(6,8)&1,7,10,13,16\\
 &(11,12)&2,6,8,15\\
 &(9,12)&0,11,12\\
 &(8,8)&3\\
 &(7,8)&5,14\\
 &(10,12)&9\\
 &(4,4)&4\\\hline
21&(10,12)&3,7,8,9,11,14,16,17,18,22\\
 &(11,12)&2,4,5,6,10\\
 &(9,12)&15,19,20,21,23\\
 &(21,24)&13
\end{array}
\]
Indeed, if $P_j=\prod_{i\in I}H_{ij}$, deleting $i$ replaces the column product $P_j$ by $P_jH_{ij}$. The same count formulas above then apply to $H^t$ with the remaining row set. 
\endgroup
The table of grouped bounds uses \eqref{eq:energy}--\eqref{eq:mod16} and the improved witnesses. For each row let $\mathcal I$ be the indicated set of column counts $1\le n\le N$, let $h=4\lceil N/4\rceil$, and set $L_n=\min\{U_0(N,n,h),U_0(n,N,h),L^{\rm imp}_n\}$, where $L^{\rm imp}_n$ is the improved bound when supplied and $+\infty$ otherwise. \begingroup
For $(N,n)=(17,9)$ take $L^{\rm imp}_9=51$ by Lemma~\ref{lem:paley51}; for $(25,21)$ take $L^{\rm imp}_{21}=117$ by Lemma~\ref{lem:paley117}.
\endgroup
 Set $E_n=B(N,n,h)$ except for $E_{10}=156$ when $N=13$. The entries $b,d$ are respectively
\begin{equation}\label{eq:groupmax}
 \max_{n\in\mathcal I}\left(\frac{E_n}{Nn}+\frac{B(n,N,h)}{N^2}\right),\qquad
 \max_{n\in\mathcal I}\frac{L_n^2}{N^2n}.
\end{equation}
The indices $n_b,n_d$ attain the respective maxima in \eqref{eq:groupmax}. For instance, when $N=17,n=9$, $B(17,9,20)=177$, $B(9,17,20)=329$ \begingroup
and $L_9=51$ by Lemma~\ref{lem:paley51}, giving $b=1990/867$ and $d=1$.
\endgroup
 The last column is the smaller of the minima of $H_3-W$ on $[0,7/10]$ and $H_1-W$ on $[7/10,1]$. These minima are evaluated at the endpoints and, for a convex quadratic $At^2+Bt+C$, also at its vertex if it lies in the interval, where its value is $C-B^2/(4A)$.
\begingroup
\begin{center}\small
\begin{tabular}{rrlllll}\toprule
$N$&$n$&$a$&$b$&$d$&$(n_b,n_d)$&margin\\\midrule
7&all&$8/7$&$110/49$&$361/343$&$(7,7)$&$236/8575$\\
10&all&$6/5$&$12/5$&$128/125$&$(10,10)$&$163/12500$\\
11&all&$12/11$&$262/121$&$1369/1331$&$(11,11)$&$1712/33275$\\
13&$\ne10$&$16/13$&$412/169$&$512/507$&$(12,6)$&$11719/1267500$\\
13&10&$16/13$&$2024/845$&$882/845$&$(10,10)$&$51/32500$\\
14&all&$8/7$&$16/7$&$576/539$&$(14,11)$&$153/13475$\\
17&$\ne9$&$20/17$&$676/289$&$2178/2023$&$(16,14)$&$9919913/12571500000$\\
17&9&$20/17$&$1990/867$&$1$&$(9,9)$&$17669/361250$\\
18&all&$10/9$&$20/9$&$784/729$&$(18,9)$&$83/18225$\\
21&all&$8/7$&$1004/441$&$475/441$&$(20,19)$&$32/11025$\\
22&all&$12/11$&$24/11$&$1156/1089$&$(22,9)$&$503/27225$\\
25&all&$28/25$&$1396/625$&$3364/3125$&$(24,20)$&$11/3125$\\
29&all&$32/29$&$1852/841$&$900/841$&$(28,16)$&$207/21025$\\
33&all&$12/11$&$2372/1089$&$3249/3025$&$(32,25)$&$18/3025$\\
37&all&$40/37$&$2956/1369$&$1458/1369$&$(36,32)$&$513/34225$\\\bottomrule
\end{tabular}
\end{center}
\endgroup

The positivity of the margins and \eqref{eq:lower}, \eqref{eq:H} give $F_A(p,q)^2<C_5(p,q)^2$ in all the listed sizes. For $5\times2$, use three rows $(1,1)$ and two rows $(1,-1)$; its Gram has diagonal 5 and off-diagonal 1, giving $a=6/5,b=56/25,d=18/25$, likewise strictly below \eqref{eq:H}. By transposition, only orders $3,5,6,9$ remain.

\begingroup
\subsection{Integer bounds used in the grouped comparisons}
The following table gives $L_n$ in increasing order of $n$. They use the two parity bounds in \eqref{eq:parity}, the improved bounds in the preceding tables, and Lemmas~\ref{lem:paley51} and \ref{lem:paley117}. Thus all the grouped $d$ bounds can be checked by squaring the displayed integers and dividing by $N^2n$. For $b$, the energy bound is obtained directly from
\[
 B(m,n,h)=b_0-\begin{cases}(b_0-m)\bmod8,&n\text{ odd},\\b_0\bmod4,&n\text{ even},\end{cases}
 \qquad b_0=\min(hn-(h-m)(n\bmod2),mn^2),
\]
where the remainder is in $\{0,\ldots,k-1\}$ for modulus $k$; use $E_{10}=156$ at $N=13$. For every indicated $n$, the comparisons to a grouped row are precisely $E_nN+B(n,N,h)n\le bN^2n$ and $L_n^2\le dN^2n$. The index pair in that row gives equality in each respective comparison.
\begin{center}\small
\begin{tabular}{rrl}\toprule
$N$&$n$&$L_n$ in increasing order of $n$\\\midrule
7&1--7&7,8,11,14,15,16,19\\
10&1--10&10,12,16,20,22,24,26,28,30,32\\
11&1--11&11,12,17,22,23,26,29,30,33,34,37\\
13&1--12&13,16,21,26,29,32,33,34,39,42,41,42\\
13&13--13&45\\
14&1--12&14,16,22,28,30,34,38,40,42,44,48,50\\
14&13--14&52,50\\
17&1--12&17,20,27,30,37,42,45,48,51,54,57,60\\
17&13--17&63,66,67,68,71\\
18&1--12&18,20,28,36,38,44,48,50,56,56,60,64\\
18&13--18&66,68,70,74,76,78\\
21&1--12&21,24,33,34,45,52,57,60,65,68,71,74\\
21&13--21&77,78,83,84,89,92,95,96,99\\
22&1--12&22,24,34,44,46,52,58,60,68,68,74,76\\
22&13--22&80,84,86,90,92,96,98,100,102,104\\
25&1--12&25,28,39,42,53,60,67,70,77,80,85,88\\
25&13--24&91,92,99,98,105,110,113,116,117,120,123,126\\
25&25--25&129\\
29&1--12&29,32,45,60,61,70,77,80,89,92,97,102\\
29&13--24&105,112,113,120,121,126,129,132,137,140,143,146\\
29&25--29&149,152,155,158,161\\
33&1--12&33,36,51,68,69,80,87,92,101,104,111,116\\
33&13--24&119,126,129,136,137,144,147,150,155,158,163,164\\
33&25--33&171,172,177,180,183,186,189,192,195\\
37&1--12&37,40,57,76,77,88,97,102,113,116,123,128\\
37&13--24&133,142,143,152,153,160,163,168,173,176,181,184\\
37&25--36&189,192,197,200,203,208,211,216,217,222,223,228\\
37&37--37&231\\\bottomrule
\end{tabular}
\end{center}
The positive margins use the two polynomials $H_3-W$ and $H_1-W$ on the closed intervals given above. For example, for $(N,n)=(17,9)$ the values at the four interval endpoints are $2593/17000$, $145109/2890000$, $17669/361250$, and $2/25$. In this case the vertex of either convex quadratic is outside its corresponding interval. The remaining rows use exactly the same endpoint-and-vertex rule.
\endgroup
\section{Orders three and six and narrow rectangles}\label{sec:small}
Write $\ell(t)=\log L(2/(1+t))^2$, $g(t)=\log R(2/(1+t))^2$. Both are convex by the perspective of log-sum-exp. The tangent to $g$ at $1/3$ gives
\begin{equation}\label{eq:g}
 g(t)>.312-.294t.
\end{equation}
Here and throughout finite decimals are exact rationals. If the four corner bounds satisfy $a\le4/3,bc\le(11/9)^2,d\le25/27$, then \eqref{eq:corners} and $\log(4/3)<.288$, $\log(11/9)<.201$, $\log(25/27)<-.076$ give $\log F_A^2\le .288-.174t-.190t^2$. Its difference from the lower line \eqref{eq:g}, with the latter first, is $(24-120t+190t^2)/1000>0$.

A single column has normalized norm one. For two columns use multiplicities $\lceil N/2\rceil,\lfloor N/2\rfloor$ of $(1,1),(1,-1)$. If $N$ is even, the columns are orthogonal. If odd, the corner squares are $(N+1)/N,(N+1)/N,(N^2+1)/N^2,(N+1)^2/(2N^2)$, by writing the transposed image as $(S+T,S-T)$ where $S,T$ are the sign sums in the two groups. For $N=3$ the bounds $a\le4/3$, $bc\le(11/9)^2$, and $d\le25/27$ hold, and for $N\ge5$ the upper polynomial in \eqref{eq:corners} is at most $6/5-4t/25-8t^2/25$, strictly below $H_1$.

For three columns use the four row types of product one, a tetrahedron with pairwise inner products $-1$ and sum of outer products $4I$. For $N=3,6,9$ take multiplicities $(1,1,1,0),(2,2,1,1),(3,2,2,2)$ respectively. \begingroup
The maximum energy and absolute sum of a direct sign image are $N+8\max c_i$ and $N+2\max c_i$: the input agrees up to sign with exactly one of the four row types. For the transposed energy, if $S_i$ denotes the sum of the signs attached to the $c_i$ copies of the $i$th type, then the squared norm is $4\sum S_i^2-(\sum S_i)^2$. This convex quadratic attains its maximum on the box $|S_i|\le c_i$ at a vertex, and every vertex is attainable. The maximum is therefore $4\sum c_i^2-\min_{\epsilon_i=\pm1}(\sum\epsilon_ic_i)^2$, where the minimum absolute sums for the three multiplicity vectors are $1,0,1$.
\endgroup Spectral squares are at most $4,8,11$, by deleting one tetrahedral row, appending two rows (whose outer-product sum has norm four), and doubling the tetrahedron and appending one row. The resulting corners are
\[
\begin{array}{c|cccc}N&a&b&c&d\\\hline
3&4/3&11/9&11/9&25/27\\
6&4/3&11/9&10/9&25/27\\
9&11/9&11/9&83/81&25/27
\end{array}
\]
In each case $a\le4/3$, $bc\le(11/9)^2$ and $d\le25/27$. Thus $\log F_A^2<g(t)$ by \eqref{eq:g} and the corner bounds.

\begingroup
For $6\times4$, choose three distinct rows of each of the two orthogonal parity bases. If $u$ and $v$ are the two omitted rows, the column Gram matrix is $8I_4-uu^t-vv^t$, where the rows are regarded here as column vectors. In particular, it is bounded above by $8I_4$. A sign input belongs, up to sign, to one of the two bases; its image in that basis has at most one coordinate of modulus four, the others being zero, and all three coordinates in the other basis have modulus two. Thus the energy is at most $16+3\cdot4=28$ and the absolute sum at most $4+3\cdot2=10$. The transposed energy is at most $8\cdot6=48$. Consequently $a\le4/3,b\le28/24=7/6,c\le48/36=4/3,d\le100/144=25/36$.
\endgroup Using $\log(14/9)<.442,\log(25/36)<-.364$, the gap in \eqref{eq:g} is at least $(24-160t+518t^2)/1000>0$.

For $6\times6$, delete two rows and two columns from $H_8$ whose intersection is Hadamard of order two. The two omitted rows restricted to the six columns are orthogonal, with inner products $S+T,S-T$ for $S,T$ sums of three signs; their energy is at least four. The transpose has the same property. Thus $a\le4/3,b,c\le44/36$, and for a six-coordinate even profile, $(z_i-2)(z_i-4)\ge0$ gives $6\sum z_i\le44+48=92$. The even absolute sum is at most 14, so $d\le49/54<25/27$.

For $6\times5$ index $H_8$ by $\mathbb F_2^3$. Keep a plane $V$ and one column $w\notin V$; omit a pair of rows differing by the nonzero member of $V^\perp$. The plane contribution is constant on each such pair, and its four values are a Walsh transform of four signs. The extra column contributes opposite signs on each pair. After omitting one pair the possible profiles are $(1,1,1,1,1,1)$, $(5,3,1,1,1,1)$, $(3,3,3,1,1,1)$. The second dominates in $r$-power sum. Consequently
\[
 J(t)=\log F_A(\infty,p_t)^2=(1+t)\log\frac{5^{2/(1+t)}+3^{2/(1+t)}+4}{6}-\log5.
\]
It is convex, with upper bounds $.237,.021,-.223$ at $0,1/3,1$. Interpolation gives $\log F_A^2\le(1-t)\log(4/3)+tJ(t)$. The upper secants of $J$ are $.237-.648t$ and $.143-.366t$ respectively. Subtraction from \eqref{eq:g} gives $(24-243t+648t^2)/1000$ and $(24-149t+366t^2)/1000$, positive by negative discriminants. This excludes order six, and all the exclusions apply to transposes by the diagonal bound.
\section{Rectangles of order nine}
In addition to \eqref{eq:g}, tangents to $g$ at $1/3$ and $1$, and to $\ell$ at $1/5$, give
\begin{equation}\label{eq:lines}
 g(t)>.3128-.2936t,\quad g(t)>.2215-.1447t,\quad \ell(t)>.3777-.5598t.
\end{equation}
If $J(y)$ bounds $\log F_A(\infty,p_y)^2$, separate convexity also gives
\begin{equation}\label{eq:Jinterp}
 \log F_A(p_x,p_y)^2\le(1-x)((1-y)\log a+y\log c)+xJ(y).
\end{equation}
For a Hadamard submatrix, the energy of a sign image is $hn$ minus the energy in the omitted rows. For even $n$, if the omitted rows have products of both signs, at least four energy units are lost by the modulo-four argument.

For $1\le r\le2$ the quadratic $P_r$ interpolating $s^r$ at $1,3,5$ majorizes $s^r$ at every positive odd integer, with nonnegative coefficients of $s^2$ and $s$. Indeed the interpolation remainder has nonpositive sign beyond $5$ since the third derivative is nonpositive. The quadratic coefficient is nonnegative by convexity, and the linear coefficient is nonnegative because $5^r+2\le3\cdot3^r$, by convexity of $(5/3)^r+2(1/3)^r$ and the endpoints $1,2$. Thus a profile of the same length, with entries in $\{1,3,5\}$ and dominating the two moments $\sum z_i,\sum z_i^2$, also dominates $\sum z_i^r$.

Use the order-twelve Paley matrix described in Section~\ref{sec:had}, with zero-based indices, and choose the following rows $I$ and columns $J$:
\begin{center}\small
\begin{tabular}{rll}\toprule
$n$&$I$&$J$\\\midrule
4&4,7,0,8,6,2,9,3,10&5,9,7,11\\
5&6,7,1,2,10,8,11,5,4&7,5,9,10,1\\
6&4,9,7,8,3,1,11,10,5&5,4,9,11,0,7\\
7&1,10,2,5,6,9,0,4,3&5,4,11,2,7,9,10\\
8&0,6,11,10,3,1,5,4,9&0,11,2,5,10,1,9,8\\\bottomrule
\end{tabular}
\end{center}
In each case $a=4/3$ and $c\le(96+n)/81$ by the omitted energy of the transposed sign image.
\begingroup
For the row-type calculations, normalize the first entry in each retained row to $1$ and encode a row $w$ by $\sum_{j=0}^{n-1}2^j\mathbf1_{w_{j+1}=-1}$. In the order of the rows $I$ above the codes are as follows. The last column gives the products of the three omitted rows before normalization, in increasing order of their indices.
\begin{center}\small
\begin{tabular}{rll}\toprule
$n$&row codes&omitted products\\\midrule
4&$8,2,0,10,14,2,6,4,8$&$-,+,-$\\
5&$18,4,26,28,8,20,22,14,24$&$+,-,-$\\
6&$24,36,22,12,34,60,48,10,58$&$+,-,+$\\
7&$52,78,98,30,118,120,0,68,18$&$-,+,-$\\
8&$0,82,106,150,160,60,216,204,116$&$-,-,-$\\\bottomrule
\end{tabular}
\end{center}
The entries follow from the quadratic residues $1,3,4,5,9$ modulo eleven and the formula for the order-twelve Paley matrix. For even $n$, the parity of the number of ones in the binary expansion of a code gives the product of that row, since row normalization does not change its product. For odd $n$, the codes still distinguish the types modulo sign.
\endgroup

For four columns the retained parity families have sizes three and six. The three types are distinct up to sign, and no type in the six occurs more than twice. The omitted rows have both parities, so $\|Au\|_2^2\le44$. The sign image has absolute sum at most $\max(12+4,6+8)=16$, by the two-basis property. Thus $b=44/36,c=100/81,d=256/324$ are bounds. Using $\log(bc)<.412,\log d<-.235$, the difference between \eqref{eq:g} and the logarithmic bound in \eqref{eq:corners} exceeds $(24-130t+359t^2)/1000>0$.

For five columns the nine row types are distinct up to sign. Each even subset of five elements has exactly five neighbors at distance four. A profile has at most one five and at most five threes, and energy at most $60-3=57$. If a five occurs, at most three threes occur. Otherwise use $2\cdot3^r\le5^r+1$. Thus the profile sum is at most $5^r+3\cdot3^r+5$. Set
\[
 J(t)=(1+t)\log\frac{5^{2/(1+t)}+3\cdot3^{2/(1+t)}+5}{9}-\log5.
\]
Its convexity and values at $0,1/3,1$ give the upper secants $.237-.492t$ and $.167-.282t$. In \eqref{eq:Jinterp}, with $x=y=t$, use $\log c<.221$ and $\log a<.288$. The differences from \eqref{eq:g} are bounded below by $(24-176t+425t^2)/1000$ and $(24-106t+215t^2)/1000$, both positive.

For six columns the omitted rows have both parities, and the retained families have sizes eight and one. Energy is at most $68$. If the family of eight has moduli zero or four, the last modulus is two or six and the sum is at most eighteen. If the family of eight has moduli two or six, a six prevents another six or a four in the last coordinate by the energy bound; otherwise the total is at most twenty. Take $b=68/54,c=102/81,d=400/486$. Since $\log(bc)<.462$ and $\log d<-.194$, the gap from \eqref{eq:g} is greater than $(24-180t+368t^2)/1000>0$.

For seven columns the energy is at most $81$. Normalize rows and input sign vector to product one, so the signed coordinates are $7,3,-1,-5$. If their multiplicities are $a_0,b_0,c_0,d_0$ then $6a_0+b_0+3d_0\le9$. If the absolute sum is at least 25, then $3a_0+b_0+2d_0\ge8$, giving $a_0=0,d_0\le1$. In the indicated row order, their sum $s$ and inner products are
\[
 s=(3,1,1,-3,1,7,1),\qquad (\langle s,A_i\rangle)_i=(7,11,7,11,11,-1,11,7,7).
\]
If $d_0=0$, then $b_0\ge8$ and $\langle s,u\rangle\ge23>\|s\|_1=17$. If $d_0=1$, then $b_0=6$ and $\langle s,u\rangle=11$. But a row with inner product $-5$ agrees with $-u$ except at one coordinate $j$, so $\langle s,u\rangle=-\langle s,A_i\rangle+2s_j(A_i)_j\le7$: for the row with inner product $-1$ with $s$, the maximum of $s_j(A_i)_j$ is three; otherwise use seven. This is impossible. Hence the absolute sum is at most $23$. The quadratic majorant $P_r$ yields the majorizing profile $z_7=(5,5,3,3,3,1,1,1,1)$.

For eight columns the energy is at most $96$ and the parity split is eight and one. If the family of eight has moduli two or six, at most two sixes occur. According to their number, the profile is coordinatewise bounded (after rearrangement) by $(6,6,2^6,0),(6,4,2^7)$ or $(8,2^8)$, where superscripts indicate multiplicities. The first dominates the other two in $r$-power sum by $6^r\ge4^r+2^r$ and $2\cdot6^r\ge8^r+2\cdot2^r$. The latter follows from the nonpositive third divided difference at $0,2,6,8$. If the family of eight has moduli zero, four or eight, the last coordinate is two or six; the possible maximal power sums are bounded by $5\cdot4^r+2^r,8^r+4^r+2^r,3\cdot4^r+6^r$. The second is bounded by the first. For the first comparison with $2\cdot6^r+6\cdot2^r$, use $2\cdot3^r-5\cdot2^r+5\ge0$: the function is convex on $[1,2]$, equals one at one, and has derivative there greater than $-1$. For the third use $3^r+6\ge3\cdot2^r$, following from $3^r+1\ge2\cdot2^r$ and $2^r\le4$ (indeed $3^r+6-3\cdot2^r\ge5-2^r>0$). Thus the majorizing profile is $z_8=(6,6,2,2,2,2,2,2,0)$.

For $n=7,8$, put $J_n(t)=(1+t)\log(\sum z_{n,i}^{2/(1+t)}/9)-\log n$. It is convex. Let $j_n$ be the piecewise affine interpolation of the following strict upper bounds and take $c_n$ as indicated:
\[
\begin{array}{c|rrrr|r}n&0&1/5&2/5&1&c_n\\\hline
7&.25132&.15802&.08323&-.06937&.24028\\
8&.28769&.17025&.07764&-.11778&.24995
\end{array}
\]
Here $c_n>\log((96+n)/81)$. By \eqref{eq:Jinterp} the logarithmic upper bounds for the two orientations are $V_n(x,y)$ and $V_n(y,x)$, where
\[
 V_n(x,y)=(1-x)((1-y).28769+yc_n)+xj_n(y).
\]
Both are nonincreasing in each variable: $j_n$ is decreasing and lies below $(1-y).28769+yc_n$ at all nodes, hence throughout. Assume $x\ge y$. For $0\le y\le1/5$, the derivative of $V_n(y,y)$ is greater than $-.258$; the first lower line for $\ell$ in \eqref{eq:lines} minus this quadratic decreases, and the gap at $1/5$ is at least $.26574-.2581636>0$. For $2/5\le y\le1$, its derivative is less than $-.315$, so the gap with $.3128-.2936y$ increases, and at $2/5$ is at least $.19536-.1946124>0$. For $1/5\le y\le2/5\le x$ it suffices to take $x=2/5$; both orientations are affine in $y$. At $y=1/5$ the maximum is at most $.2361852<.3128-.2936/5$, and at $y=2/5$ the previous gap applies.

\begingroup
To detail the rational comparisons, write $j_n(t)=\alpha+\beta t$ on a given interpolation interval and put $a=.28769$, $c=c_n$. Then
\[
 V_n(t,t)=a+(c-2a+\alpha)t+(a-c+\beta)t^2.
\]
The derivative is affine, with the following endpoint bounds (the listed decimals enclose the exact rational values obtained by substitution):
\begin{center}\small
\begin{tabular}{c|rr|rr|rr}\toprule
&\multicolumn{2}{c|}{$[0,1/5]$}&\multicolumn{2}{c|}{$[1/5,2/5]$}&\multicolumn{2}{c}{$[2/5,1]$}\\
$n$&left&right&left&right&left&right\\\midrule
7&-.08378&-.251416&-.232906&-.363522&$<-.31567$&$<-.56398$\\
8&-.03774&-.257524&-.232694&-.402818&-.347878&-.69343\\\bottomrule
\end{tabular}
\end{center}
In the middle triangle the partial derivatives in $x$ of $V_n(x,y)$ and $V_n(y,x)$ are respectively $j_n(y)-a-(c-a)y$ and $(1-y)(c-a)+y\beta$. Both are decreasing in $y$ there. Their values at $y=1/5$ are respectively $-.120188,-.112718$ for $n=7$ and $-.109892,-.122802$ for $n=8$. This verifies the uniform derivative estimates used in the comparison.
\endgroup
It remains to compare on $1/5\le y\le x\le2/5$. For $E_1$ use the fixed vectors
\[
 v=(-20,-12,9,-9,-16),\qquad w=(-20,-17,-14,-20,-14),
\]
with $w^tE_1v=3170$. Since $E_1$ is equivalent to its transpose, $\log F_{E_1}(p_x,p_y)^2\ge2\log3170-\log125-f_v(y)-f_w(x)$, where $f_z(s)=(1-s)\log(\sum |z_i|^{2/(1-s)}/5)$ is convex. Taking the secants over $[1/5,2/5]$ yields the lower affine bound
\[
 P(x,y)=.37596-.35677y-.09523x.
\]
The class $E_2$ has profile $(3,3,3,3,1)$, whose normalized logarithmic lower bound throughout this triangle exceeds $.34$ (evaluate at $y=2/5$), which exceeds $P$. By Proposition~\ref{prop:square}, $P$ is therefore a lower bound for $\log C^2$. In the triangle the derivatives in $x$ of both orientations of $V_n$ are at most $-.1098$, hence the gap with $P$ increases in $x$. Along the diagonal the derivative of $V_n$ is greater than $-.403$, so the gap with $P$ decreases, attaining a lower bound $.37596-.452(2/5)-.1946124=.0005476>0$ at the right endpoint. This completes the rectangular exclusions.

\section{The square of order nine}
\subsection{The construction and its sign profiles}
Let $G=(\mathbb Z/3\mathbb Z)^2$, $S=\{(0,0),(1,0),(0,1)\}$, and $Q=(1-2\mathbf1_{i-j\in S})_{i,j\in G}$. Let $K=\{(0,0),(1,1),(2,2)\}$ and $K_1,K_2,K_3$ its cosets. The six nonzero differences of elements of $S$ are exactly the elements outside $K$. Thus two distinct translates intersect once if the indices are in different cosets, and are disjoint otherwise. Consequently $12I-Q^tQ=W^tW$, with columns of $W$ equal within each coset and taken from the three types $(1,1,1),(1,-1,-1),(-1,1,-1)$. If $T_i=\sum_{j\in K_i}v_j$ then
\begin{equation}\label{eq:gram9}
 \|Wv\|_2^2=4\sum_iT_i^2-(\sum_iT_i)^2\ge\tfrac13(\sum_jv_j)^2.
\end{equation}
Also $Q^t$ is equivalent to $Q$ by inversion of the indices.
\begingroup
Indeed, for $i\ne j$ the Gram entry is $9-4|S|+4|(i+S)\cap(j+S)|$, namely $-3$ in the same coset and $1$ otherwise; its diagonal is $9$. The proposed columns of $W$ have squared norm $3$ and mutual inner product $-1$, proving the matrix identity entry by entry. In particular $Q^tQ\le12I$. The formula \eqref{eq:gram9} follows from the Gram matrix $4I_3-J_3$ and Cauchy--Schwarz.
\endgroup
\begin{lemma}\label{lem:two}
For a sign vector $\sigma$, $d=Q^t\sigma$ has odd coordinates, $\sum|d_i|\le27$, $\sum d_i^2\le105$. If both bounds are attained, its signed profile, up to global sign and permutation, is $(-3,-3,-3,1,1,1,5,5,5)$, and this profile occurs.
\end{lemma}
\begingroup
\begin{proof}
Since each coordinate of $W\sigma$ is the sum of nine signs, its modulus is at least one. The Gram identity gives $\|Q\sigma\|_2^2=108-\|W\sigma\|_2^2\le105$. The same holds for $Q^t$ by equivalence. All coordinates are odd.

It suffices to consider $Q\sigma$, with at most four negative signs, since a global sign change and inversion of the indices preserve the assertion. Write $T\subset G$ for the negative positions, $k=|T|$, and $a_i=\#\{j\in T:i\in j+S\}$. Then
\[
 (Q\sigma)_i=3-2k+4a_i,\qquad 0\le a_i\le3,\qquad \sum_i a_i=3k.
\]
Put $n_j=\#\{i:a_i=j\}$ for $0\le j\le3$. The identities
\[
 \sum_{j=0}^3n_j=9,\quad \sum_{j=0}^3 jn_j=3k,\quad n_2+3n_3=\sum_{\{u,v\}\subset T}|(u+S)\cap(v+S)|
\]
follow by counting points, incidences, and unordered pairs of incidences. The final sum counts pairs in different cosets of $K$.

For $k=0$, the image is constant three and has energy $81$. For $k=1$ the absolute sum is $3\cdot5+6=21$. For $k=2$ the two translates are disjoint or meet once, yielding absolute sums $21$ or $23$. None of these cases attains both bounds.

For $k=3$, at most three pairs meet. The union has at least six points, since $9-(n_1+n_2+n_3)=n_2+2n_3\le3$. The absolute sum is $9+6n_0\le27$. Equality forces $n_0=3$ and $n_2+2n_3=3$. Since $n_2+3n_3\le3$, we have $n_3=0$ and $n_2=n_1=3$. The signed coordinates are therefore three each of $-3,1,5$, with energy $105$. This case occurs for $T=S$: the translates indexed by $(0,0),(1,0),(0,1)$ meet pairwise at three different points, and have no common point.

For $k=4$, at least two indices lie in one coset; hence at most five pairs meet. Since $n_2+2n_3\le n_2+3n_3\le5$, the union has at least $12-5=7$ points, so $n_0\le2$. The absolute sum is $3+10n_0+2n_1$ (use the values $-5,-1,3,7$ and the two incidence identities). If $n_0=2$, the equations imply $n_2+2n_3=5$, so $n_3=0$, $n_2=5$, $n_1=2$. The absolute sum is $27$, but the energy is $2\cdot25+2+5\cdot9=97$. If $n_0=1$, then $n_1=4+n_3\le6$, and the absolute sum is at most $25$. If $n_0=0$, the absolute sum is at most $3+2\cdot9=21$. 
\end{proof}
\endgroup

\begin{lemma}\label{lem:quadratic}
If $f$ is increasing, convex, of class $C^3$ on $[0,\infty)$, with $f'''\le0$ and $f'(0)\ge0$, the interpolating quadratic at $1,3,5$ majorizes $f$ at all positive odd integers and has nonnegative quadratic and linear coefficients.
\end{lemma}
\begin{proof}The interpolation remainder is nonpositive for $s>5$. With $e_1=(f(3)-f(1))/2,e_2=(f(5)-f(3))/2$, the coefficients are $(e_2-e_1)/4$ and $2e_1-e_2$. For the latter use that $f'(s)/s$ is nonincreasing, from concavity of $f'$ and $f'(0)\ge0$. Integration on $[1,3]$ and $[3,5]$, comparing with $sf'(3)/3$, gives $e_2\le2e_1$.\end{proof}
By the moment bounds and the quadratic comparison for powers (or regularization of $s^r$ in Lemma~\ref{lem:quadratic}),
\begin{equation}\label{eq:Qedge}
 \|Q\|_{\infty\to r}^r=3(5^r+3^r+1),\quad 1\le r\le2.
\end{equation}

\subsection{Gram compensation and two interior estimates}
Fix $2<p<\infty$, $r=p'$, $z_0>0$, $a=(r-1)z_0^{r-2}$, $b=(2-r)z_0^{r-1}$ and $\lambda>0$. Weighted AM--GM gives $z^r\le az^2+bz$. Put
\[
 \beta(d)=\max_{x\ge0}(12ax^2+dx-\lambda x^p),\qquad d\ge0.
\]
The unique positive maximizer $x_d$ solves $\lambda p x_d^{p-1}=24ax_d+d$ (the left side after division by $x_d$ minus $24a+d/x_d$ is strictly increasing). At $d=0$ the positive maximizer is also strictly better than zero. \begingroup
Implicit differentiation gives $\beta'(d)=x_d$ and, with $D_d=\lambda p(p-1)x_d^{p-2}-24a>0$,
\[
 \beta''(d)=D_d^{-1},\qquad
 \beta'''(d)=-\frac{\lambda p(p-1)(p-2)x_d^{p-3}}{D_d^3}<0.
\]
Indeed $\lambda p x_d^{p-2}=24a+d/x_d$, so $D_d\ge24a(p-2)>0$. These formulas extend to $d=0$ from the right.
\endgroup Thus Lemma~\ref{lem:quadratic} applies to $f(s)=\beta(bs)$.
Write $B_j=\beta(jb)$ for $j=1,3,5$ and
\begin{align*}
 U_1&=9\lambda+5B_1+4B_5,\\
 U_2&=9\lambda+2B_1+5B_3+2B_5,\\
 U_3&=9\lambda+3\{\beta(|b-a|)+\beta(3b+a)+\beta(|5b-a|)\}+3a/4.
\end{align*}
For $\sum|v_i|^p\le9$, take $\sigma$ such that $\|Qv\|_1=\langle Q^t\sigma,v\rangle$. With $d=Q^t\sigma$, the scalar inequality and \eqref{eq:gram9} yield
\[
 \|Qv\|_r^r\le12a\|v\|_2^2-a\|Wv\|_2^2+b\langle d,v\rangle.
\]
\begingroup
The energy is congruent to $9$ modulo eight, and the absolute sum is odd. Thus, outside the exceptional case of Lemma~\ref{lem:two}, at least one of the two moment bounds decreases by $8$ or $2$, respectively.
\endgroup
\begingroup
In the nonexceptional case, discard $-a\|Wv\|_2^2$ and use $\sum_i|v_i|^p\le9$ to bound the right-hand side by $9\lambda+\sum_i\beta(b|d_i|)$. The absolute sum is at most $25$ or the energy at most $97$.
\endgroup Write $P(s)=As^2+Bs+D$ for the interpolator of $s\mapsto\beta(bs)$, where $A,B\ge0$. The exceptional profile has quadratic sum $3(B_1+B_3+B_5)=105A+27B+9D$. The nonexceptional profiles give at most $105A+25B+9D$ or $97A+27B+9D$. Since $8A=B_5-2B_3+B_1$ and $2B=3B_3-2B_1-B_5$, these yield $U_1$ and $U_2$, respectively, after adding $9\lambda$. \begingroup
In the exceptional case replace $(v,\sigma)$ by $(-v,-\sigma)$ if necessary, so that $d$ has three coordinates of each of the values $-3,1,5$. All terms in the inequality are unchanged. Since \eqref{eq:gram9} bounds the Gram term by the square of the total sum of the coordinates of $v$, their positions in the three cosets are immaterial.
\endgroup
 With $T=\sum v_i$, use $-aT^2/3+aT\le3a/4$ in \eqref{eq:gram9} to replace the linear term by $\langle bd-a\mathbf1,v\rangle+3a/4$. Applying the definition of $\beta$ coordinatewise gives $U_3$. Therefore for $t=1-2/p$,
\begin{equation}\label{eq:H9}
 H(t):=\log F_Q(p_t,p_t)^2\le (1+t)\log(\max_i U_i/9)-\log9.
\end{equation}

\subsection{Scalar comparisons}
For $s>0$, convexity of $u^{p/2}$ gives $x^p\ge (p/2)s^{p-2}x^2-(p/2-1)s^p$. Completing the square, provided the denominator is positive, yields
\begin{equation}\label{eq:beta}
 \beta(d)\le\lambda(p/2-1)s^p+\frac{d^2}{2\lambda ps^{p-2}-48a}.
\end{equation}
Choose $z_0=19/5$ and $(p,\lambda)=(5/2,509/100),(10/3,59/25)$. Direct rational substitutions in \eqref{eq:beta}, specified in Section~\ref{sec:arithmetic}, give
\[
\begin{array}{c|rrrr}t&U_1<&U_2<&U_3<&H(t)<\\\hline
1/5&67.8163&68.6640&68.1625&.2413\\
2/5&48.3132&49.5451&49.4478&.1908
\end{array}
\]
Also $H(0)\le\log(4/3)<.28769$, $H(1)=0$ by \eqref{eq:Qedge}. The same formula gives
\[
 J(t)=\log F_Q(\infty,p_t)^2=(1+t)\log\frac{5^{2/(1+t)}+3^{2/(1+t)}+1}{3}-\log9,
\]
with upper bounds $.25952,.1861,.1261$ at $0,1/5,2/5$. For any $0\le u<v\le1$, let $h_u,h_v,j_u$ be upper bounds for $H(u),H(v),J(u)$. Separate convexity and symmetry give at $(u,v)$ the bound $c=h_u+(v-u)(j_u-h_u)/(1-u)$, and for $t\in[u,v]$, $w=(t-u)/(v-u)$,
\[
 H(t)\le(1-w)^2h_u+2w(1-w)c+w^2h_v.
\]
On the intervals $[0,1/5],[1/5,2/5],[2/5,1]$ this yields respectively
\[
 .28769-.05634t-.87805t^2,\quad .246+.091t-.5725t^2,\quad(2248-713t-1535t^2)/9000.
\]
The tangent to $g$ at $1/2$ gives $g(t)>.2913-.2416t$. The first polynomial is strictly below the lower line for $\ell$ in \eqref{eq:lines} (the decreasing gap at $1/5$ is $.02444$); the second is below $.3128-.2936t$ since the gap $.0668-.3846t+.5725t^2$ has negative discriminant. The third is below $.2913-.2416t$ since the gap times 90000 is $3737-14614t+15350t^2$, again with negative discriminant. By monotonicity in both variables, $F_Q(p,q)\le F_Q(p_t,p_t)$ for $t=\min(1-2/p,1-2/q)$, so the square of order nine is strictly below \eqref{eq:lower}.

The maximum in \eqref{eq:finite-reduction} is attained only at the sizes $5\times5$, $5\times4$, $5\times3$ and their transposes. The fixed-size minimizations in Section~\ref{sec:classification} give \eqref{eq:band}, proving Theorem~\ref{thm:continuous}.
\section{Boundary formulas and regions with a single representative}
\subsection{The Hilbertian boundary}
For $2<p<\infty$, maximizing $Q_2$ on $B_p^5$ is equivalent to maximizing $Q_2(x)-\sum_i|x_i|^p$ on $\mathbb R^5$: on the ray $x=tv$, $\|v\|_p=1$, the maximum is
\[
 (1-2/p)(2/p)^{2/(p-2)} Q_2(v)^{p/(p-2)}.
\]
Writing $(\sum x_i)^2=\max_u(2u\sum x_i-u^2)$, for the maximizing value of $u$ the two pairs may be chosen equal since they maximize the same separate function. Hence a maximizer of the form $(a,b,a,b,z)$ exists. The identity $Q_2(a,b,a,b,z)=Q_1(z,-b,a,a,-b)$ and equality of the $\ell_p$ norms prove that $E_1$ minimizes the square norm on this boundary. Continuity gives the endpoints $p=2,\infty$.

Put $\theta=1-2/p$. A sign vector attaining the energy 37 and the eigenvector $(1,-1,1,-1,0)$ give
\begin{equation}\label{eq:hilblow}
 F_E(p,2)^2\ge\max\{37/25,(8/5)(4/5)^\theta\}.
\end{equation}
The Hadamard estimate in Section~\ref{sec:had} bounds all $N\ne5$ by $4/3<37/25$ in squared normalized norm. For $5\times n$, $n=1,2,3$, the spectral squares of the witnesses are $5,6,7$ (the last has Gram $4I_3+J_3$ up to sign changes), so both orientations are below $37/25$. For $4\times5$ the corners for $D_1^t$ give $F_{D_1^t}(p,2)^2\le(8/5)(7/10)^\theta\le(8/5)(4/5)^\theta$.
For $5\times4$, the Gram $G=D_2^tD_2$ satisfies
\[
 G=\begin{pmatrix}5&1&1&3\\1&5&1&-1\\1&1&5&-1\\3&-1&-1&5\end{pmatrix}.
\]
Direct expansion gives
\[
 3x^t(\operatorname{diag}(8,20/3,20/3,8)-G)x
 =(3x_1-3x_4-x_2-x_3)^2+4(x_2-x_3)^2.
\]
For $p>2$, with $s=p/(p-2)=1/\theta$, H\"older gives
\[
 F_{D_2}(p,2)^2\le\frac85\left(\frac{1+(5/6)^s}{2}\right)^{1/s}.
\]
For $s\ge3$, the mean inside parentheses is at most $341/432<4/5$. For $1\le s\le3$, the power mean is at most $(341/432)^{1/3}$ and $(8/5)(341/432)^{1/3}<37/25$ since $341/432<(37/40)^3$. Now apply \eqref{eq:hilblow}, using continuity at $p=2$. This proves \eqref{eq:edge2} of Theorem~\ref{thm:boundary}, including $C(2,2)=\sqrt{8/5}$ and $C(\infty,2)=\sqrt{37}/5$. Monotonicity proves optimality of the uniform constant.

\subsection{The boundary with an infinite exponent}
The square minimum is $L(r)$ by the profiles in Section~\ref{sec:classification}. The $5\times3$ minimum is $R(r)$ because all four row types occur, with one repeated. \begingroup
For $B^t$, the first four signs of an input give three coordinates of a Walsh transform of order four; the fifth sign is added to each of them. A concentrated transform gives a profile bounded by $(5,1,1)$, and a flat transform by $(3,3,3)$. The latter is attained, for instance at $(1,1,1,-1,1)^t$. Since $25+1+1=3\cdot9$ and the vector $(25,1,1)$ majorizes $(9,9,9)$, concavity of $s^{r/2}$ gives $5^r+2\le3\cdot3^r$. Hence $F_{B^t}(\infty,r')=R(1)\le R(r)$.
\endgroup

For $D_2$, there are two rows of one parity and three of the other, all distinct modulo sign. Thus $(4,2,2,2,0)$ is attainable and dominates every profile. It gives
\[
 F_{D_2}(\infty,r')=\left(\frac{2^r+3}{5}\right)^{1/r}\le R(r).
\]
Indeed the inequality is equivalent to $f(r)=2\cdot3^{r/2}+3\cdot3^{-r/2}-2^r-3\ge0$. Here $f(2)=0$, $f'(2)=(5/2)\log3-4\log2<0$ since $3^5<2^8$, and $f''>0$: the ratio $(2\cdot3^{r/2}+3\cdot3^{-r/2})/2^r$ decreases and is at least $7/4$, and $(7/16)(\log3)^2>(\log2)^2$ since $3^7>2^{11}$ and $11/7>4/\sqrt7$. Thus $f>0$ for $r<2$.
\begingroup
For $D_1^t$, the Walsh transform of four signs is concentrated or flat. In the latter case its four signs have negative product, and the repeated row adds a common sign to its four coordinates. The resulting absolute values are bounded in power sum by $(3,3,3,1)$, which is attained. A concentrated transform gives a profile bounded by $(5,1,1,1)$. Both displayed profiles have energy $28$, and $(25,1,1,1)$ majorizes $(9,9,9,1)$. Concavity of $s^{r/2}$ shows that $(3,3,3,1)$ gives the maximum for $1\le r\le2$.
\endgroup In particular $F_{D_1^t}(\infty,\infty)=1$, $F_{D_1^t}(\infty,2)^2=28/25$. If $1\le r\le3/2$, interpolation in $1/r$ gives $F_{D_1^t}(\infty,r')\le(28/25)^{1/3}<\sqrt{27/25}\le R(r)$ (cube the squared inequality). If $3/2\le r\le2$, the bound $\sqrt{28}/5<R(3/2)\le R(r)$ applies. Formula \eqref{eq:band} now gives \eqref{eq:edgeinfty}, including $C(\infty,\infty)=3\sqrt3/5$.

\begingroup
\begin{corollary}\label{cor:transition}
There is a unique $r_0\in(1,2)$ satisfying $L(r_0)=R(r_0)$. On the boundary with an infinite exponent, $R(r)$ gives the constant for $1\le r<r_0$, and $L(r)$ gives it for $r_0<r\le2$; at $r_0$ both sizes attain the constant. Numerically $r_0\approx1.61315$, corresponding to $p=r_0/(r_0-1)\approx2.63092$.
\end{corollary}
\begin{proof}
The sign of $L(r)^r-R(r)^r$ is the sign of
\[
 h(r)=(5/3)^{r/2}+(3/5)^{r/2}+3\cdot15^{-r/2}-2-3\cdot3^{-r}.
\]
Put $a=\tfrac12\log(5/3)>0$. Then
\[
 h'(r)=a(e^{ar}-e^{-ar})+3\{(\log3)3^{-r}-(\log\sqrt{15})15^{-r/2}\}.
\]
For $r\ge1$ we have $(\sqrt{15}/3)^r\ge e^a>1+a>1+a/\log3=\log\sqrt{15}/\log3$, since $\log3>1$. Thus both terms in $h'(r)$ are positive. At $r=1$, $L(1)=11/(5\sqrt5)<9/(5\sqrt3)=R(1)$, while $L(2)^2=37/25>7/5=R(2)^2$. Continuity and strict monotonicity yield the result.
\end{proof}
\endgroup
\subsection{Regions with a specified minimizing matrix}
The maximizing profiles for the norms $\|A\|_{\infty\to r}$, $1\le r\le2$, and spectral square bounds $h$ are as follows (for $D_1$ take the larger of the two power sums).
\[
\begin{array}{c|l|c}
 A&\text{profiles}&h\\\hline
 E&(5,3,1,1,1)&8\\
 B&(3,3,1,1,1)&7\\
 B^t&(3,3,3)&7\\
 D_2&(4,2,2,2,0)&8\\
 D_1&(4,4,0,0,0),(2,2,2,2,2)&8\\
 D_1^t&(3,3,3,1)&8
\end{array}
\]
The profiles follow from the classification in Section~\ref{sec:classification}; the spectral bound for $D_2$ follows from its Gram being four times a sum of two orthogonal projections (onto the subspaces spanned by each parity family). For a profile $z$ of an $m\times n$ matrix of maximal dimension five, write
\[
 J_z(t)=(1+t)\log(\sum_i z_i^{2/(1+t)}/m)+\log(m/(5n)).
\]
Let $J_A$ be the maximum over its listed profiles. If $p\ge2d$, monotonicity and separate convexity give
\[
 \log F_A(p,p_t)^2\le V_{A,d}(t)=d^{-1}\log(h/5)+(1-d^{-1})J_A(t).
\]
Each $V$ is convex. Thus to bound it, or a finite maximum of such functions, by an affine lower bound on an interval it suffices to check the endpoints. The profile $(3,3,3,3,1)$ gives a lower bound $u(t)>.3920-.1364t$ for the log-square of every square matrix of order five not equivalent to $E$. Indeed $E_2$ has that profile, and the other classes have the lower power sum in Lemma~\ref{lem:graph}, which is at least $4\cdot3^r+1$ because $5^r+1\ge2\cdot3^r$ for $r\ge1$. In the following table let $\mathcal A=\{B,B^t,D_2,D_1^t\}$ and $\mathcal B=\{E,B^t,D_1,D_1^t\}$. The last two columns are upper bounds for the maximum of $V_{A,d}$ over the indicated family, at the two endpoints.
\begin{center}\small
\begin{tabular}{rllrrr}\toprule
$d$&interval&family&lower line&left&right\\\midrule
12&$[1/19,1/7]$&$E$&$.3920-.1364t$&.3656&.3134\\
12&$[1/19,1/7]$&$\mathcal A$&$.3777-.5598t$&.3263&.2920\\
25&$[47/153,1/3]$&$\mathcal B$&$.3128-.2936t$&.2200&.2078\\
50&$[1/3,3/5]$&$\mathcal B$&$.3128-.2936t$&.2023&.0936\\
200&$[3/5,2/3]$&$\mathcal B$&$.3128-.2936t$&.0879&.0784\\\bottomrule
\end{tabular}
\end{center}
All endpoint comparisons are strict. In the first region the first line proves that $E$ minimizes the square norm, and the second line that all rectangular minima are strictly below it. In the other regions, $B$ minimizes the $5\times3$ norm and the listed competitors bound the minima in all remaining dimensions strictly below it. This proves Proposition~\ref{prop:regions}, including the transposed identities.

\section{Details of the scalar inequalities}\label{sec:arithmetic}
For $1\le u\le2$, with $z=(u-1)/(u+1)$,
\begin{equation}\label{eq:log}
 2\sum_{j=0}^{19}\frac{z^{2j+1}}{2j+1}\le\log u\le
 2\sum_{j=0}^{19}\frac{z^{2j+1}}{2j+1}+\frac{2z^{41}}{41(1-z^2)}.
\end{equation}
For other positive rational arguments first reduce by powers of two (or take reciprocals). Fractional powers are bounded by rational intervals using integer powers. The following table gives the lower integer endpoint of $10^5\log x$ with upper endpoint one larger:
\[
\begin{array}{c|rrrrrr}
x&2&3&5&13/5&12/5&4/3\\\hline
 &69314&109861&160943&95551&87546&28768\\
x&11/9&27/25&14/9&36/25&19/15&5/4\\\hline
 &20067&7696&44183&36464&23638&22314
\end{array}
\]
This also proves the two estimates for $k$ in Lemma~\ref{lem:tangent}. \begingroup
For \eqref{eq:H}, at $t=1$ use $\log(125/81)>2(z+z^3/3)>13/30$, $z=22/103$. Consequently $f'(1)<-(9/25)(13/30)=-.156$.
\endgroup
 At $t=1/2$, put $c=3^{1/3}$ and $M=(6c+3)/5$. Then $f(1/2)=M^{3/2}/3$ and $f'(1/2)=f(1/2)(\log M-8c\log3/(6c+3))$. Cubing and squaring give $1.44224<c<1.44225$, $\textcolor{black}{1.1860<f(1/2)<1.1862}$, and \eqref{eq:log} gives $.84616<\log M<.84618$, yielding the derivative bounds used there.

For each fixed profile the derivative of $\phi(t)=(1+t)\log(\sum z_i^{2/(1+t)}/m)$ is
\[
 \phi'(t)=\log(\sum z_i^r/m)-r\frac{\sum z_i^r\log z_i}{\sum z_i^r},\qquad r=2/(1+t).
\]
The exact tangent intercepts and slopes for $\ell$ at $1/5$, $g$ at $1/3$, and $g$ at $1/2$ lie respectively in
\[
 (.37783,.37784)\times(-.55970,-.55969),\quad
 (.31288,.31290)\times(-.29354,-.29352),
\]
\[
 (.29141,.29143)\times(-.24156,-.24155).
\]
For example at $1/3$ use $173205/10^5<\sqrt3<173206/10^5$ and $B=(6\sqrt3+3)/5$, $w=6\sqrt3/(6\sqrt3+3)$: $g(1/3)=(4/3)\log B-\log3$, $g'(1/3)=\log B-(3/2)w\log3$. The enclosures $77599/10^5<w<77600/10^5$ and $98524/10^5<\log B<98525/10^5$ yield the stated weaker lines. The tangent at one follows from the derivative formula and the log bounds; the intercept and slope for the profile $(3,3,3,3,1)$ at zero lie in $(.39204,.39205)$ and $(-.13637,-.13635)$.

For the power comparisons used in the body, each integer $M$ in the next table specifies $M\le10^8 z^r<M+1$ (check by raising both sides to the denominator of $r$). A dash indicates an unused entry.
\begin{center}\small
\begin{tabular}{r|rrrrr}
$r$&$z=2$&$z=3$&$z=4$&$z=5$&$z=6$\\\midrule
$5/3$&317480210&624025146&--&1462008869&1981156349\\
$10/7$&269180038&480398665&--&996617657&1293137313\\
$20/17$&--&364182542&--&664228723&--\\
$19/10$&373213196&806362613&1392880901&2128349806&--\\
$7/4$&336358566&683852117&1131370849&1671850762&--\\
$153/100$&288785839&537026256&833972608&1173340622&--\\
$3/2$&282842712&519615242&800000000&1118033988&--\\
$4/3$&--&432674871&--&854987973&--\\
$5/4$&237841423&394822203&565685424&747674390&--\\
$6/5$&229739670&373719281&527803164&689864830&--
\end{tabular}
\end{center}
For the interior test on $E_1$, use $v/20,w/20$ in place of $v,w$; the corresponding enclosures for $10^8(k/20)^p$ are
\[
\begin{array}{c|rrrrrr}
p\backslash k&9&12&14&16&17&20\\\hline
5/2&13584112&27885480&40996341&57243340&66611208&100000000\\
10/3&6982994&18218145&30455107&47529869&58174110&100000000
\end{array}
\]
Substitution in the secants gives an intercept in $(.375963,.375965)$ and the absolute values of the coefficients of $y,x$ in $(.356761,.356764)$ and $(.095228,.095230)$, respectively. The bilinear constant after scaling the vectors is $2\log(317/40)-\log125$. For example, the final diagonal comparison for $9\times7$ and $9\times8$ has gap at least $5476/10^7$. The endpoint comparisons for the regions with a single representative are obtained by substituting the power bounds in the formula for $J_z$ and then in $V_{A,d}$; \eqref{eq:log} bounds each logarithm by rational numbers. The polynomial comparisons require only rational subtraction.

For \eqref{eq:beta} let $u=(19/5)^{1/3}$, $k=2$, $a=2/(3u)$, $b=u^2/3$ at $p=5/2$; and $u=(19/5)^{1/7}$, $k=3$, $a=3/(7u^4)$, $b=4u^3/7$ at $p=10/3$. The inequalities $1560490750<10^9u<1560490751$ and $1210113840<10^9u<1210113841$, respectively, follow by integer powers. Choose $s=(j/10^4)^k$ in \eqref{eq:beta}, as in the table.
\[
\begin{array}{c|rr|rr}
 &\multicolumn{2}{c|}{p=5/2}&\multicolumn{2}{c}{p=10/3}\\
d&j&\beta(d)<&j&\beta(d)<\\\hline
b&8869&1.0174&9358&1.2251\\
3b&9979&2.4719&9992&3.0762\\
5b&10795&4.2298&10420&5.2369\\
b-a&8478&.6956&9273&1.0635\\
3b+a&10214&2.9075&10041&3.2771\\
5b-a&10598&3.7409&10383&5.0120
\end{array}
\]
Let $a_-,a_+,b_-,b_+$ be the bounds obtained from the indicated interval for $u$, and let $d_+$ be respectively $b_+,3b_+,5b_+,b_+-a_-,3b_++a_+,5b_+-a_-$. For each row, with $D=2\lambda p(j/10^4)^{k(p-2)}-48a_+$ and $T$ the listed upper bound, direct substitution gives $D>0$ and
\[
 [T-\lambda(p/2-1)(j/10^4)^{kp}]D>d_+^2.
\]
The powers $k(p-2)$ and $kp$ are integers, so these are rational inequalities. They imply the bounds on $U_i$ and \eqref{eq:H9} upon substitution and application of \eqref{eq:log}.
\section*{Funding}
D. M. Pellegrino acknowledges support from CNPq through Grants No. 406457/2023-9 (CNPq/MCTI No. 10/2023), 403964/2024-5 (MCTI/CNPq No. 16/2024), and 305807/2025-0.


\section*{Data availability}
No datasets were generated or analysed in this study.

\section*{Competing interests}
The competing interests declaration is to be confirmed by all authors before submission.

\section*{AI Declaration}
ChatGPT (OpenAI) was used to explore the finite-dimensional reduction, identify candidate matrices and estimates, develop and check parts of the analytical arguments, and assist in drafting and revising the manuscript.

\end{document}